\documentclass{amsart}
\usepackage{graphicx}

\newtheorem{theorem}{Theorem}[section]
\newtheorem{lemma}[theorem]{Lemma}
\newtheorem{corollary}[theorem]{Corollary}

\theoremstyle{definition}

\theoremstyle{remark}

\numberwithin{equation}{section}

\newcommand{\boundary}{\operatorname{Bd}}
\newcommand{\closure}{\operatorname{Cl}}

\newcommand{\glb}{\operatorname{glb}}
\newcommand{\horiline}{\operatorname{Horiz}}
\newcommand{\identity}{\operatorname{Id}}

\newcommand{\lub}{\operatorname{lub}}
\newcommand{\mesh}{\operatorname{mesh}}
\newcommand{\neighborhood}{\operatorname{N}}
\newcommand{\rwidth}{\operatorname{Width}}
\newcommand{\rheight}{\operatorname{Height}}
\newcommand{\sstar}{\operatorname{st}}
\newcommand{\successor}{\operatorname{Nxt}}
\newcommand{\vertline}{\operatorname{Vert}}

\newcommand{\E}{\operatorname{\mathbb E}}

\newcommand{\Z}{\operatorname{\mathbb Z}}

\begin{document}

\title[]{More Monotone Open Homogeneous Locally Connected Plane Continua}

\author[]{Carl R. Seaquist}
\address{Department of Mathematics and Statistics, Texas Tech University,
Lubbock, TX 79409-1042}
\email{seaqucr\@math.ttu.edu}
\thanks{This research was supported in part by Texas ARP grant \#003644-012.}

\subjclass{primary 54B15, secondary 54C10, 54E45, 54F50}

\begin{abstract}
This paper constructs a continuous decomposition of the Sierpi\'nski
curve into acyclic continua one of which is an arc.  This
decomposition is then used to construct another continuous decomposition
of the Sierpi\'nski curve.
The resulting
decomposition space
is homeomorphic to the continuum obtained from taking
the Sierpi\'nski curve and identifying two points on the boundary
of one of its complementary domains.
This outcome is shown to imply that there are continuum many topologically
different one dimensional locally connected plane
continua that are homogeneous with respect to monotone open maps.
\end{abstract}

\maketitle
\section{Introduction}
In a recent paper~\cite{Prajs1998}, J. R. Prajs
remarks without proof that there are infinitely many topologically different
locally connected plane continua that are homogeneous with
respect to monotone open mappings.
These spaces are described as {\em two} dimensional and
monotone open equivalent to the Sierpi\'nski curve. 
Here, using the construction techniques developed in~\cite{Seaquist1998}
to prove that the Sierpi\'nski curve is monotone open homogeneous,
we prove that there are continuum
many one dimensional locally connected plane continua that are monotone open
homogeneous.
A continuum, $X$, is {\em monotone open
homogeneous} if for any two points, $x$ and $y$, in $X$
there is a monotone open map from $X$ onto $X$ so that
$f(x)=y$. 
A {\em monotone map} is one with connected fibers and
an {\em open map} is one that preserves open sets.
By {\em map} we mean a continuous function.
We say two continua, $X$ and $Y$, are {\em monotone open equivalent}
if there is a monotone open map from $X$ onto $Y$ and vice versa.

This paper makes use of the fact that there is a continuous
decomposition of the Sierpi\'nski curve into acyclic continua,
one of which is an arc lying in the boundary of the unbounded complementary domain,
such that the decomposition space is homeomorphic to the Sierpi\'nski curve.
The details of the construction of such a decomposition are given in Section 3--5
of this paper and are
almost identical to the
continuous decomposition
described in~\cite{Seaquist1998}; the only difference being that there
in the final stage of the construction the decomposition
elements are wrapped around
a bounded complementary domain,
while here the decomposition elements are stretched
and bent in order to lie alongside an arc
on the boundary of the unbounded complementary domain.
This decomposition can be used to show that there
is a continuous decomposition of the Sierpi\'nski curve into acyclic continua
so that the decomposition space is homeomorphic to the
continuum which results from taking the Sierpi\'nski curve and
identifying two points on the boundary of a bounded complementary domain.
Combining this result with those of \cite{Seaquist1998},
we get that this last continuum is monotone open equivalent
to the Sierpi\'nski curve and so is also monotone open
homogeneous.
A generalization of this result shows that there are continuum many
locally connected plane continua that are monotone open homogeneous.
This result, however, leaves open the question asked by J. Prajs
on whether or not there might exist a nondegenerate
locally connected plane continuum that is monotone open homogeneous
but is not either a simple closed curve or monotone open
equivalent to the Sierpi\'nski curve.

\section{Main Result}
In this section we will assume the following theorem, which will be
proved in the remaining sections of the paper.
The proof of the theorem in Sections~3--5 in no way depends on the results
in Section~2.

\begin{theorem}\label{thm:construction}
There exists a continuous decomposition of the Sierpi\'nski
curve into acyclic continua, one of which is an arc lying in
the boundary of the unbounded complementary domain, such that the decomposition
space is homeomorphic to the Sierpi\'nski curve.
\end{theorem}

We will call a bounded complementary domain of the Sierpi\'nski curve
a {\em hole}.
We will call the continuum that results from taking a Sierpi\'nski
curve and identifying $n$ points on
the boundary of one hole of the curve
a {\em Sierpi\'nski curve with one pinched hole with
$n$ lobes}.
Notice that a Sierpi\'nski curve with one pinched hole has
one local cut point, which we will call the {\em center} of the pinched hole.
Given a Sierpi\'nski curve with one pinched hole
we call the union of the complementary domains which have boundaries
that intersect at the local cut point {\em a pinched hole}, while
each such complementary domain is called a {\em lobe} of the pinched hole.
We remark that using the same techniques
that are used in \cite{Whyburn1958a} by G.~T.~Whyburn
to prove that two $S$-curves are homeomorphic, we can prove that if we have two
Sierpi\'nski curves each with one pinched hole with $n$ lobes,
then they are homeomorphic.
We will also need the notion of a pinched hole with a {\em pinched lobe}.
By a Sierpi\'nski curve with one pinched hole with $n$ lobes one of which
is pinched we mean the continuum that results from identifying
two points of the boundary of a lobe of a pinched hole
where neither of the points are the center of the
pinched hole.

By the standard square Sierpi\'nski curve, we mean the Sierpi\'nski curve
that results from removing open square holes from the unit square
$[0,1]\times [0,1]$ in the standard way.
Specifically, let $S_1=I\times I$ and $S_2$ be
the continuum that results from dividing $S_1$ into $9$ identical
squares and removing the interior of the center one.
For $i\in\Z^+$, we obtain
$S_{i+1}$ by taking each of the $8^{(i-1)}$ squares that make up $S_i$ and
dividing it into $9$ identical squares and then removing the interior of the center one.
The standard square Sierpi\'nski curve is $S=\cap^\infty_{i=1}S_i$.

Using Theorem~\ref{thm:construction} we are able to prove the following theorem:

\begin{theorem}\label{thm:StoB}
Let $B$ be a Sierpi\'nski curve with one pinched hole with two lobes.
Then there is a continuous decomposition of the Sierpi\'nski curve into
acyclic continua so that the decomposition space is homeomorphic
to $B$.
\end{theorem}
\begin{proof}
Let $S$ be the standard square Sierpi\'nski curve.
From Theorem~\ref{thm:construction}
there is a continuous decomposition $G_1$ of $S$ into acyclic continua one
of which is an arc that without loss of generality we can assume is the
left edge of $S$.
We can also assume that $S/G_1$  is the Sierpi\'nski curve $T(S)$ where
$T:S\rightarrow{\E}^2:(x,y)\mapsto (x,yx)$.
In addition we would like to assume that
the members of $G_1$ that lie along the right edge of $S$ are single points.  
To see that we can make this assumption without loss of generality consider the
following argument.
Let ${g:S\rightarrow S/G_1}$ be the natural map and consider the
decomposition space $Z$ that results from collapsing to points all the members
of $G_1$ that intersect the right edge of $S$.
We will denote by $\phi$ this quotient map.
Since the right edge of $S$ is a closed set and $g$ is a closed map it can be shown
that $\phi$ is closed.
Also since $g$, in addition to being closed, is monotone
and open, the map $\pi':Z\rightarrow S/G_1$ defined by $\pi'(z)=f(\phi^{-1}(z))$ is
closed, monotone, and open.
Finally, since $\phi$ is a homeomorphism on the boundary points of the complementary domains
of $S$ it can be shown using results of~\cite{Whyburn1958a} that $Z$ is homeomorphic
to the Sierpi\'nski curve.
Now $\pi'$ induces a continuous decomposition of $Z$ into acyclic
continua so that in addition to one element of the decomposition lying along the boundary
of the unbounded complementary domain there is another arc of the unbounded complementary
domain that is the union of degenerate elements of the decomposition.
Thus we can assume without loss of generality
that the members of $G_1$ that lie along the right edge of $S$ are single points.
Let $\pi_1$ be the quotient map from $S$ to $S/G_1$.

Let
$G_2$ be a continuous decomposition of $S$ into acyclic continua with the quotient map $\pi_2$ so that
$S/G_2=T'(S)$ where 
$T':S\rightarrow{\E}^2:{(x,y)\mapsto (x,y(1-x))}$; so that $\pi_2^{-1}((1,0))$
is the right edge of $S$; and so that $\pi_2$ is the identity on the
left side of $S$.  See Figure~\ref{fig:StoB}.

\begin{figure}[tb]
\begin{center}
\includegraphics*{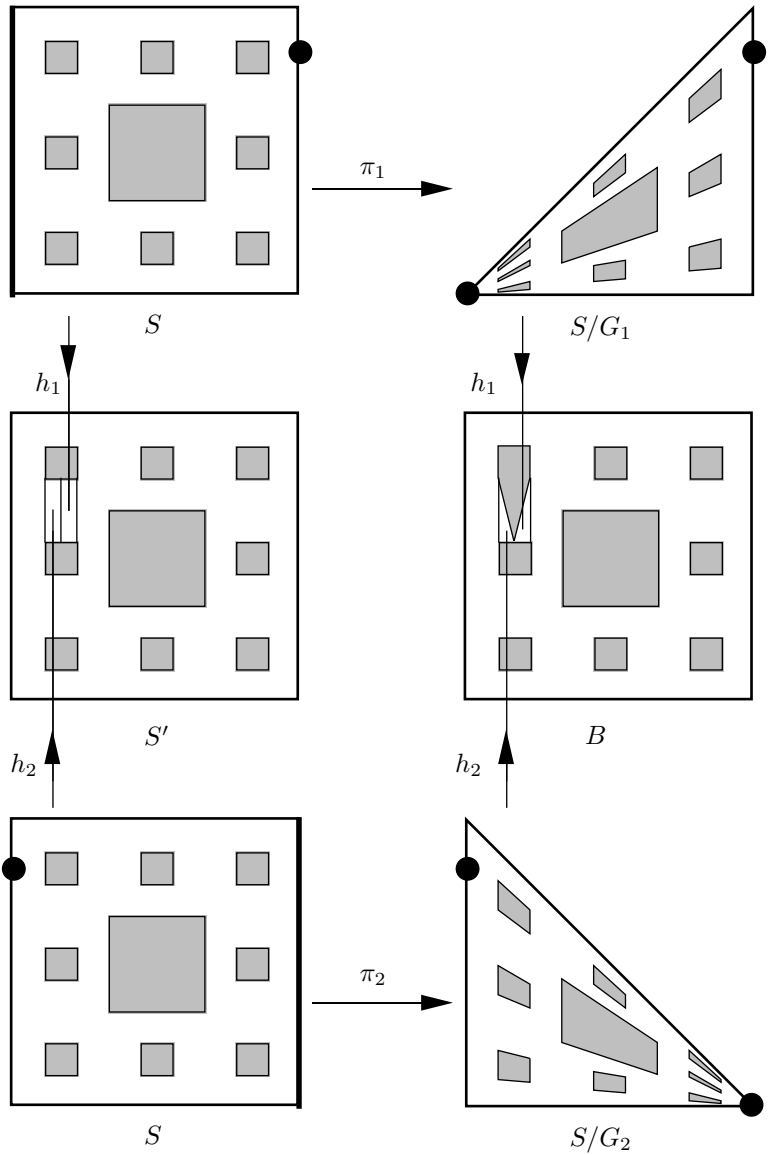}
\end{center}
\setlength{\unitlength}{.01in}
\hspace*{-4.25in}\begin{picture}(0,0)
\put(200,525){\makebox(0,0)[tl]{$\pi_1$}}
\put(87,445){\makebox(0,0)[tl]{$S$}}
\put(310,445){\makebox(0,0)[tl]{$S/G_1$}}
\put(30,415){\makebox(0,0)[tl]{$h_1$}}
\put(258,415){\makebox(0,0)[tl]{$h_1$}}
\put(87,230){\makebox(0,0)[tl]{$S'$}}
\put(318,230){\makebox(0,0)[tl]{$B$}}
\put(17,215){\makebox(0,0)[tl]{$h_2$}}
\put(250,215){\makebox(0,0)[tl]{$h_2$}}
\put(200,102){\makebox(0,0)[tl]{$\pi_2$}}
\put(87,20){\makebox(0,0)[tl]{$S$}}
\put(310,20){\makebox(0,0)[tl]{$S/G_2$}}
\end{picture}
\caption{There is a monotone open map from $S$ to $B$.}
\label{fig:StoB}
\end{figure} 
Let $h_1:S\rightarrow ([\frac{1}{6},\frac{2}{9}]\times [\frac{5}{9},\frac{7}{9}]):
(x,y)\mapsto(\frac{1}{18}x+\frac{1}{6},\frac{2}{9}y+\frac{5}{9})$
and $h_2:S\rightarrow ([\frac{1}{9},\frac{1}{6}]\times [\frac{5}{9},\frac{7}{9}]):
(x,y)\mapsto(\frac{1}{18}x+\frac{1}{9},\frac{2}{9}y+\frac{5}{9})$.
Denote by $S'$ the Sierpi\'nski curve that results from removing the
rectangles $([\frac{1}{9},\frac{1}{6}]\times [\frac{5}{9},\frac{7}{9}])$ and
$([\frac{1}{6},\frac{2}{9}]\times [\frac{5}{9},\frac{7}{9}])$ from $S$
and replacing them with
$h_2(S)$ and $h_1(S)$ respectively.
Denote by $B$ the Sierpi\'nski curve with one pinched hole with two lobes that
results from removing the  
rectangles $([\frac{1}{9},\frac{1}{6}]\times [\frac{5}{9},\frac{7}{9}])$ and
$([\frac{1}{6},\frac{2}{9}]\times [\frac{5}{9},\frac{7}{9}])$
from $S$ and replacing them with
$h_2(S/G_2 )$ and $h_1(S/G_1)$ respectively.
We can now define the map $f:{\E}^2\rightarrow {\E}^2$ as follows:
\[ f(x,y) = \left\{
\begin{array}{ll}
h_1\circ \pi_1\circ h_1^{-1}(x,y)&\mbox{if $(x,y)\in([\frac{1}{6},\frac{2}{9}]\times [\frac{5}{9},
\frac{7}{9}])$;} \\
h_2\circ \pi_2\circ h_2^{-1}(x,y)&\mbox{if $(x,y)\in([\frac{1}{9},\frac{1}{6}]\times [\frac{5}{9},
\frac{7}{9}])$;} \\
(x,y)&\mbox{otherwise.}
\end{array}
\right. \]
Thus $f|S'$ is a continuous function from $S'$ onto $B$
which is monotone and open because $\pi_1$ and
$\pi_2$ are monotone open maps.
Therefore the theorem holds.
\end{proof}

\begin{corollary}\label{cor:StoBn}
Let $B$ be a Sierpi\'nski curve with one pinched hole with $n$ lobes with $n\in\Z^+$.
Then there is a continuous decomposition of the Sierpi\'nski curve into
acyclic continua so that the decomposition space is homeomorphic
to $B$.
\end{corollary}

The proof this corollary is very similar in style to the proof of
Theorem~\ref{thm:StoB};
however, here we must repeat the construction described
$n-1$ more times.
In each of these constructions we collapse an arc that has one end point at the center
of a pinched hole and another end point on the boundary of a different
hole and avoids all other boundary points of holes.

\begin{corollary}\label{cor:StoBnp}
Let $B$ be a Sierpi\'nski curve with one pinched hole with $n$ lobes
one of which is pinched  where $n\in\Z^+$.
Then there is a continuous decomposition of the Sierpi\'nski curve into
acyclic continua so that the decomposition space is homeomorphic
to $B$.
\end{corollary}

To prove this corollary we can use the techniques used in the proof of Theorem~\ref{thm:StoB}
to collapse an arc that has one end point on the boundary of a lobe (other than the center point of a pinched
hole) and the other end point on the boundary of
another hole and avoids all other boundary points of holes. 

\begin{theorem}\label{cor:moequivalence}
If $B$ is a Sierpi\'nski curve with a pinched hole with two lobes, then
there is a monotone open map from $B$ onto a Sierpi\'nski curve.
\end{theorem}
\begin{proof}
From the construction described in~\cite{Seaquist1998}
we know that there is a monotone open map, $g$, from the Sierpi\'nski curve,
$S$, onto $S$ so that for some $p\in S$, $g^{-1}(p)$ is the boundary of a hole.
Consider the identification $\phi$ of two points of $g^{-1}(p)$.
Now $\phi$ is a closed map and it can be shown using techniques used
in~\cite{Whyburn1958a} that the image of $S$ under $\phi$ is homeomorphic to $B$.
Let $h$ be homeomorphism from $B$ onto $\phi(S)$.
Since $g$ is open and closed the map $g\circ \phi^{-1}\circ h$ is a monotone open map
from $B$ to $S$.
\end{proof}

We now have the following two corollaries which can be proved as above.

\begin{corollary}\label{cor:BntoS}
If $B$ is a Sierpi\'nski curve with a pinched hole with $n$ lobes, then
there is a monotone open map from $B$ onto a Sierpi\'nski curve.
\end{corollary}

\begin{corollary}\label{cor:BnptoS}
If $B$ is a Sierpi\'nski curve with a pinched hole with $n$ lobes one of which is pinched, then
there is a monotone open map from $B$ onto a Sierpi\'nski curve.
\end{corollary}

\begin{figure}[tb]
\begin{center}
\includegraphics*{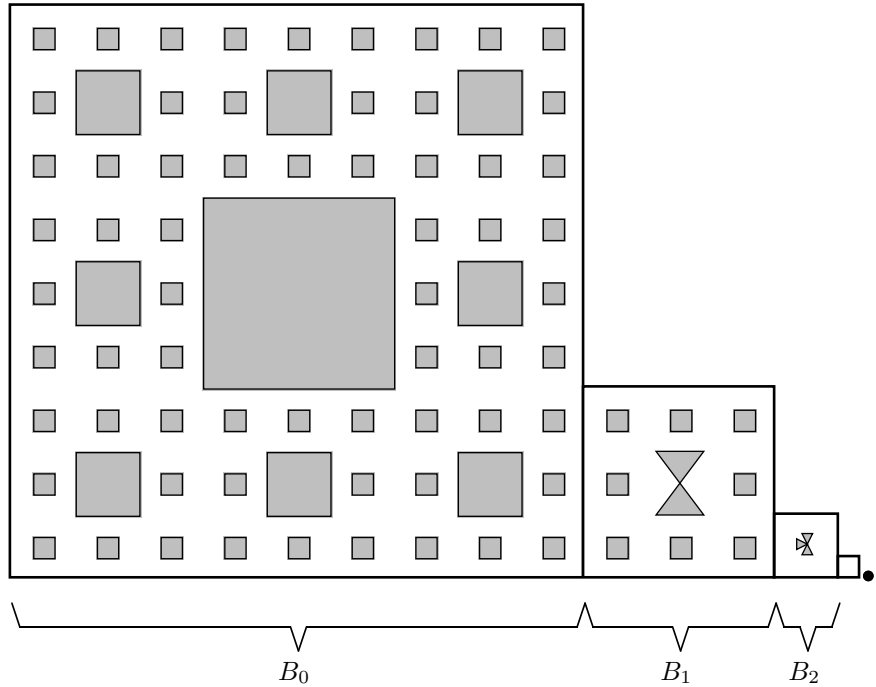}
\end{center}
\setlength{\unitlength}{.01in}
\hspace*{-4.7in}\begin{picture}(0,0)
\put(150,15){\makebox(0,0)[tl]{$B_0$}}
\put(350,15){\makebox(0,0)[tl]{$B_1$}}
\put(417,15){\makebox(0,0)[tl]{$B_2$}}
\end{picture}
\caption{Figures like $B=\cup_{i=0}^\infty B_i \cup \{(\frac{3}{2},0)\}$
are monotone open homogeneous.}
\label{fig:Bunion}
\end{figure}

Now consider a Sierpi\'nski curve with a sequence of pinched holes.
Specifically, let $B_0$ be the standard square Sierpi\'nski curve
and for each $n\in\Z^+$ let $B_n$ be a Sierpi\'nski
curve with one pinched hole with $n+1$ lobes so that the boundary
of the unbounded component of the complement is the same as the boundary
of the rectangle
$$\Bigl[\sum^{n-1}_{i=0}\Bigl(\frac{1}{3}\Bigr)^i,
\sum^n_{i=0}\Bigl(\frac{1}{3}\Bigr)^i\Bigr]\times
\Bigl[0,\Bigl(\frac{1}{3}\Bigr)^n\Bigr].$$
See Figure~\ref{fig:Bunion}.
Let $B$ be the closure of $\cup B_n$.
We will now show that $B$ is monotone open homogeneous.
Let $S_0$ be the standard Sierpi\'nski curve and $S_n$ be a Sierpi\'nski
curve with the same unbounded complementary domain as $B_n$.  Let $S$ be the closure of $\cup S_n$.
By Whyburn's characterization of the Sierpi\'nski curve~\cite{Whyburn1958a} we know that $S$
is a Sierpi\'nski curve.  From Corollary~\ref{cor:StoBn} there is a monotone open map
$f_n$ from $S_n$ onto $B_n$ for each non-negative integer.
In fact there exists such a map $f_n$ that is the identity on the boundary of
the unbounded complementary domain of $S_n$.
Define $f:S\rightarrow B$ as follows: 
\[ f(x) = \left\{
\begin{array}{ll}
f_n(x)&\mbox{if $x\in S_n$;} \\
x&\mbox{otherwise.}
\end{array}
\right. \]
Now $f$ is continuous, monotone, and open.
From Corollary~\ref{cor:BntoS} we can similarly construct a monotone open map $g$ from $B$
onto $S$.  Thus $B$ is monotone open equivalent to $S$.  Since $S$ is
monotone open homogeneous~\cite{Seaquist1998} we have that $B$ is monotone open homogeneous.

Now consider $r\in[0,1]$ and let $0.b_1b_2b_3...$ be the binary representation of $r$.
We can create a new continuum $B^r$ by either pinching or not pinching one of
the lobes of the pinched hole of $B_n$ depending on whether or not $b_n$ is $0$
or $1$.
For each $r\in[0,1]$ using corollaries~\ref{cor:StoBn}, \ref{cor:StoBnp}, \ref{cor:BntoS},
and~\ref{cor:BnptoS} we can show
that $B^r$ is monotone open equivalent to $S$ and so is monotone open homogeneous.
Since $B^r$ is not homeomorphic to $B^{r'}$ unless $r=r'$, we have proved the following
theorem.

\begin{theorem}
There are continuum many topologically different locally connected plane continua
that are monotone open homogeneous.
\end{theorem}

The remainder of the paper describes a construction that proves
Theorem~\ref{thm:construction}.  This construction follows
that described in~\cite{Seaquist1998} very closely and
in fact, much is identical but
is included here to simplify the reading of this paper.
Note that J.~Prajs has suggested~\cite{Prajs1999} an alternate proof to Theorem~\ref{thm:construction}
based on the fact that there is a continuous decomposition of the
Sierpi\'nski curve into pseudo-arcs and on a lemma about extending functions
on Sierpi\'nski curves.

\section{Preliminaries and Overview of Construction}

Before giving an overview of our construction we introduce the following
notation.
If $P$ is a collection of sets,
then $P^*$ denotes the union of members of
$P$.
If $p$ is a set, then $\sstar^1(p,P)=\{p'\in P:p'\cap p \ne\emptyset\}$ and
inductively $\sstar^{i+1}(p,P)=\sstar^{1}(\sstar^{i}(p,P)^*,P)$. We
abbreviate $\sstar^1(p,P)$ by $\sstar(p,P)$.
If $U\subset {\E}^2$, then by $U^c$ we mean the complement
of $U$ with respect to ${\E}^2$.
We will use $\vertline(x)$ to denote
the vertical line running through the point
$(x,0)$ and $\horiline(y)$ to denote a horizontal line
running through the point $(0,y)$.
If $U$ is a bounded subset of ${\E}^2$ then we define
$\rwidth(U)$ to be $\lub\{x|(x,y)\in U\}-\glb\{x|(x,y)\in U\}$.
Similarly we define
$\rheight(U)$ to be
$\lub\{y|(x,y)\in U\}-\glb\{y|(x,y)\in U\}$.

We will make use of the construction lemma,
Lemma~7 in~\cite{Seaquist1995b}, which is originally due to Lewis and
Walsh \cite{Lewis1978}.
We will construct a sequence $\{Y_n\}_{n=1}^\infty$ of
continua so that $Y_{n+1}\subset Y_n$ and so that
$Y=\cap_{n=1}^\infty Y_n$
is the Sierpi\'nski curve.
We will set $Y_0={\E}^2$ and $Y_1=D$ where
$D$ is the rectangle $[\frac{1}{2},1]\times [0,2]$.
We pick this particular rectangle to be compatible with~\cite{Seaquist1998}.
We will denote by $E$ the left edge of $D$.
While constructing $Y_n$, we will describe
a partition $P_n$ of $Y_n$ into cells 
with nonoverlapping interiors along with a bending function
$f_n$ on $D$.
The cells $P_n$ will be constructed to satisfy the constraints
imposed by Lemma~7 in~\cite{Seaquist1995b} so that
$G=\{\cap_{n=1}^\infty \sstar(p_n,P_n)^*
: \cap_{n=1}^\infty p_n \ne \emptyset
\mbox{ where } p_n\in P_n\}$
is a continuous decomposition
of $Y$.
All the members of $G$ that do not intersect
$E$ will be nondegenerate
and those that do intersect $E$ will be single points.
From the functions
$\{f_n\}_{n=1}^\infty$
we will construct a homeomorphism $F:(Y_1\backslash E)\rightarrow
(Y_1\backslash E)$ that bends
the elements of $G$ back and forth close to the left edge of $D$
thus creating a decomposition of $F(Y\backslash E)\cup E$;
namely,
$$G'=\{F(g)|g\in G \mbox{ and } g\cap E=\emptyset\}\cup\{E\}.$$
The functions
$\{f_n\}_{n=1}^\infty$
will have been carefully defined so that the holes in $Y$ will not be stretched
too much and so $F(Y\backslash E)\cup E$ will be the Sierpi\'nski curve.
The decomposition $G'$ under the quotient topology will be shown to be
homeomorphic to the Sierpi\'nski curve.

Rather than describe the cells of $P_n$ directly,
we first describe a partition $Q_n$ of
a continuum $X_n$ into fairly simple
cells with nonoverlapping interiors
and then take $Y_n$ to be
the image of $X_n$ under a homeomorphism $H_n$.
The collection
$P_n$ is obtained by applying the same homeomorphism $H_n$ to each of the
cells
$q_n\in Q_n$.
In the construction, $X_n$ will be a continuum
with finitely many complementary regions with disjoint boundaries that
are simple closed curves.

The construction starts with $X_1=Y_1=D$ and proceeds
inductively.
Note that $E\subset\vertline(1/2)$.
Assuming we are at stage~$n$ of the construction,
we are given the continuum
$X_n\subset D$
and
$\widehat{R}_n$, which is
either a horizontal or vertical division of $D$.
A vertical (resp. horizontal) division of $D$ is the collection
$R=\{[\frac{1}{2}+(i-1)a,\frac{1}{2}+ia]\times[0,2]:i\in\{1,...,1/(2a)\}\}$ where
$1/(2a)\in\Z^+$
(resp. $R=\{[\frac{1}{2},1]\times[(i-1)a,ia]:i\in\{1,...,2/a\}\}$ where
$2/a\in\Z^+$).  The mesh of $R$ is $a$ and each member of $R$
is called a strip of $R$. 

We will call the bounded complementary regions of $X_n$ {\it holes}.
We will be careful to insure that the
outside edges of $X_n$ will coincide with the outside edges of $D$.
We will also make sure that
the boundary of $X_n$ is a finite number of disjoint simple closed
curves and that the holes of $X_n$ are open squares with edges
parallel to
either the $x$-axis or the $y$-axis.
For the continua $X_n$ that arise in the construction
we use {\em vertical boundary} to mean the left and right edges of $D$ unioned
with the vertical line segments that make up the left and right edges of the
holes of $X_n$.
The term {\em horizontal boundary} is used in a similarly manner.

We now give an overview of the construction at stage~$n$.
Given a vertical (respectively horizontal)
division $\widehat{R}_n$ we refine it to obtain $R_n$.
The common part of each strip of $R_n$ and $X_n$ is then
partitioned into cells
with nonoverlapping interiors to obtain the collection of
cells $Q_n$.
Once $Q_n$ is defined, a homeomorphism $h_n:D\rightarrow D$ is defined
so that
$\{h_n^{-1}(q_n):q_n\in Q_n\}$ is a collection of identical rectangles with
nonoverlapping interiors whose union is $X_n$.
We will define $h_n$ so that it will leave the boundary of $X_n$ invariant
and will be the identity on the left edge of $D$.
We then set $Y_n=h_1\circ\cdot\cdot\cdot\circ h_{n-1}(X_n)$ and the
collection
$P_n$ is defined to be
$\{h_1\circ\cdot\cdot\cdot\circ h_{n-1}(q_n):q_n\in Q_n\}$.
The homeomorphism $h_1\circ\cdot\cdot\cdot\circ h_{n-1}$
will be denoted by $H_n$.
To continue to stage $n+1$ we use
$\{h_n^{-1}(q_n):q_n\in Q_n\}$ to define
$\widehat{R}_{n+1}$, a horizontal (respectively vertical) division of $D$.
Thus the construction alternates between working with horizontal
and vertical divisions of $D$.
Arbitrarily, we let $\widehat{R}_n$ be a vertical division when
$n$ is odd and a horizontal division when
$n$ is even.
To continue to stage $n+1$ we must also define $X_{n+1}$.
For
some of the $q_n\in Q_n$
(exactly which ones will be made clear later)
we define
a small open $s_n$ by $s_n$ square hole,
referred to as $w_n$, that will be centered in
the rectangle $h_n^{-1}(q_n)$.
The parameter $s_n$ is a rational number which helps control the
construction
at stage~$n$.
We set
$
W_n=\{w_n:\exists q_n \in Q_n
$ and $
w_n
$ is an open $s_n$ by $s_n$ square centered in $
h_n^{-1}(q_n)\}
$.
We will define $X_{n+1}$ to be
$X_{n}\backslash W^*_{n}$.
Thus $X_{n+1} = D\backslash (\cup_{i=1}^n
W^*_{i})$.

\begin{figure}[tb]
\begin{center}
\includegraphics*{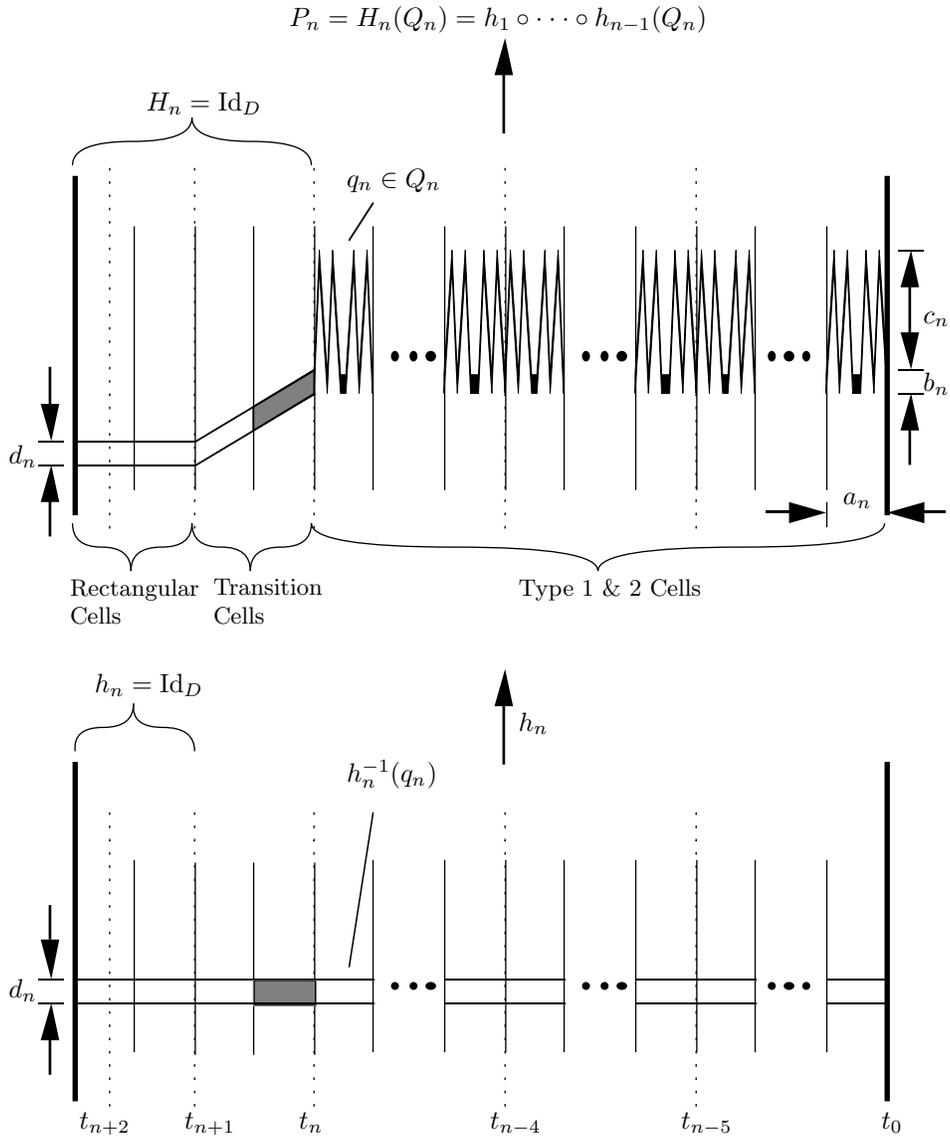}
\end{center}
\setlength{\unitlength}{.01in}
\hspace*{-5.15in}\begin{picture}(0,0)
\put(50,37){\makebox(0,0)[tl]{$t_{n+2}$}}
\put(105,37){\makebox(0,0)[tl]{$t_{n+1}$}}
\put(165,37){\makebox(0,0)[tl]{$t_{n}$}}
\put(265,37){\makebox(0,0)[tl]{$t_{n-4}$}}
\put(365,37){\makebox(0,0)[tl]{$t_{n-5}$}}
\put(470,37){\makebox(0,0)[tl]{$t_{0}$}}
\put(58,267){\makebox(0,0)[tl]{$h_n=\identity_D$}}
\put(280,247){\makebox(0,0)[tl]{$h_n$}}
\put(45,317){\makebox(0,0)[tl]{\small Rectangular}}
\put(45,302){\makebox(0,0)[tl]{\small Cells}}
\put(120,317){\makebox(0,0)[tl]{\small Transition}}
\put(120,302){\makebox(0,0)[tl]{\small Cells}}
\put(280,317){\makebox(0,0)[tl]{\small Type 1 \& 2 Cells}}
\put(85,572){\makebox(0,0)[tl]{$H_n=\identity_D$}}
\put(160,617){\makebox(0,0)[tl]
{$P_n=H_n(Q_n)=h_1\circ\cdot\cdot\cdot\circ h_{n-1}(Q_n)$}}
\put(13,107){\makebox(0,0)[tl]{$d_n$}}
\put(13,387){\makebox(0,0)[tl]{$d_n$}}
\put(450,362){\makebox(0,0)[tl]{$a_n$}}
\put(492,457){\makebox(0,0)[tl]{$c_n$}}
\put(492,425){\makebox(0,0)[tl]{$b_n$}}
\put(190,532){\makebox(0,0)[tl]{$q_n\in Q_n$}}
\put(190,222){\makebox(0,0)[tl]{$h_n^{-1}(q_n)$}}
\end{picture}
\caption{Cell types in a
vertical division.  The homeomorphism, $h^{-1}_n$, straightens
the upper and lower boundaries of the cells.}
\label{fig:verti}
\end{figure}

\begin{figure}[tb]
\begin{center}
\includegraphics*{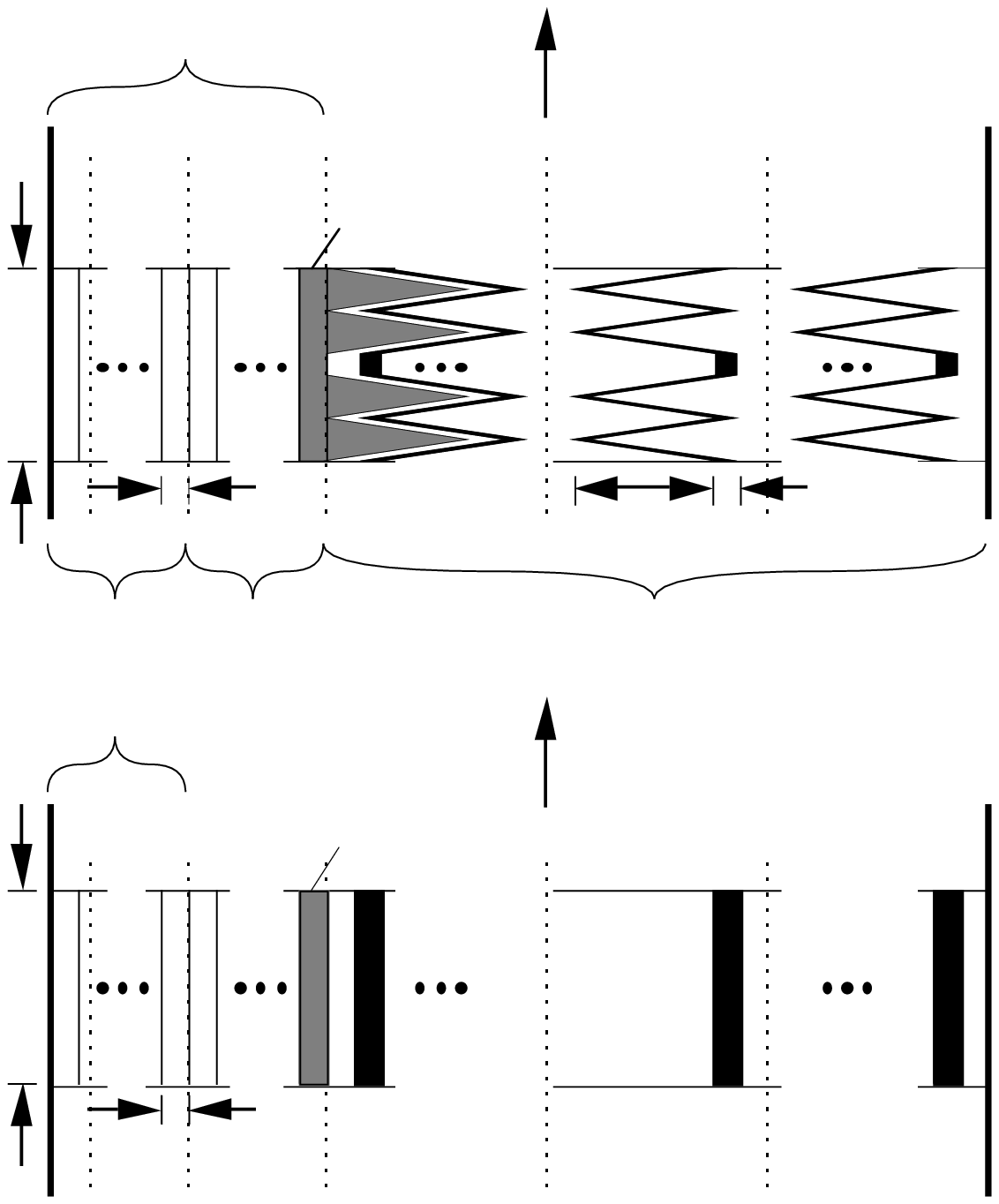}
\end{center}
\setlength{\unitlength}{.01in}
\hspace*{-4.75in}\begin{picture}(0,0)
\put(30,37){\makebox(0,0)[tl]{$t_{n+2}$}}
\put(85,37){\makebox(0,0)[tl]{$t_{n+1}$}}
\put(145,37){\makebox(0,0)[tl]{$t_{n}$}}
\put(245,37){\makebox(0,0)[tl]{$t_{n-4}$}}
\put(345,37){\makebox(0,0)[tl]{$t_{n-5}$}}
\put(450,37){\makebox(0,0)[tl]{$t_{0}$}}
\put(38,264){\makebox(0,0)[tl]{$h_n=\identity_D$}}
\put(270,247){\makebox(0,0)[tl]{$h_n$}}
\put(25,307){\makebox(0,0)[tl]{\small Rectangular}}
\put(25,292){\makebox(0,0)[tl]{\small Cells}}
\put(100,307){\makebox(0,0)[tl]{\small Transition}}
\put(100,292){\makebox(0,0)[tl]{\small Cells}}
\put(260,307){\makebox(0,0)[tl]{\small Type 1 \& 2 Cells}}
\put(65,572){\makebox(0,0)[tl]{$H_n=\identity_D$}}
\put(160,597){\makebox(0,0)[tl]
{$P_n=H_n(Q_n)=h_1\circ\cdot\cdot\cdot\circ h_{n-1}(Q_n)$}}
\put(10,142){\makebox(0,0)[tl]{$a_n$}}
\put(10,422){\makebox(0,0)[tl]{$a_n$}}
\put(80,72){\makebox(0,0)[tl]{$d_n$}}
\put(82,352){\makebox(0,0)[tl]{$d_n$}}
\put(295,352){\makebox(0,0)[tl]{$c_n$}}
\put(337,352){\makebox(0,0)[tl]{$b_n$}}
\put(170,507){\makebox(0,0)[tl]{$q_n\in Q_n$}}
\put(170,222){\makebox(0,0)[tl]{$h_n^{-1}(q_n)$}}
\end{picture}
\caption{Various types of cells also occur
in a horizontal division. The right most transition cell
is shaded.  No Type~2 cells are shown.}
\label{fig:horiz}
\end{figure}

Because we want the members of $G$ to have smaller and smaller
diameters towards the left of edge of $D$,
the left most cells $q_n\in Q_n$ will be rectangles;
whereas the right most cells will be as described in \cite{Seaquist1995b}.
More specifically, we define the sequence $t_n=1/2+(1/2)^{n+1}$ for
$n\in \{0,1,2,...\}$ and constrain the construction of cells $Q_n$ at stage
$n$ so that all cells to the left of $\vertline(t_{n+1})$ are
simply rectangles and
those to the right of $\vertline(t_n)$ are
of the two types described in \cite{Seaquist1995b}.
Cells between $\vertline(t_{n+1})$ and
$\vertline(t_n)$ are called transition
cells.
Additionally, the construction
is constrained at stage $n$ by five positive rational
constants, $a_n$, $a'_n$, $b_n$, $c_n$, and $s_n$
and a positive integer
$k_n$.
To facilitate our discussion of the application of these parameters we
give an informal description of the cells $q_n\in Q_n$.
At each stage there are four kinds of cells $q_n$ in $Q_n$:
rectangular, transition,
Type~1, and Type~2.
Rectangular cells occur to the left of
$\vertline(t_{n+1})$ and when $n$ is odd are
$a_n$ by $d_n$ rectangles; i.e., they have width
$a_n$ and height $d_n$.  When $n$ is even they are $d_n$ by $a_n$ rectangles.
Transition cells in $Q_n$ occur between $\vertline(t_{n+1})$ and
$\vertline(t_{n})$ and are different in nature
depending on whether we are building
vertical cells ($n$ odd) or horizontal cells ($n$ even).
If we are building vertical cells,
then all the transition cells are trapezoids with their
left and right  vertical boundaries being parallel.
These cells have the potential to be
relatively tall.
See Figure~\ref{fig:verti}.
If we are building horizontal cells,
then all except the right most transition cells
are $d_n$ by $a_n$ rectangles.
The right most transition cells have a left boundary
that is a straight line and a right boundary
that is the left boundary of a Type~1 cell.
See Figure~\ref{fig:horiz}.

\begin{figure}[tb]
\begin{center}
\includegraphics*{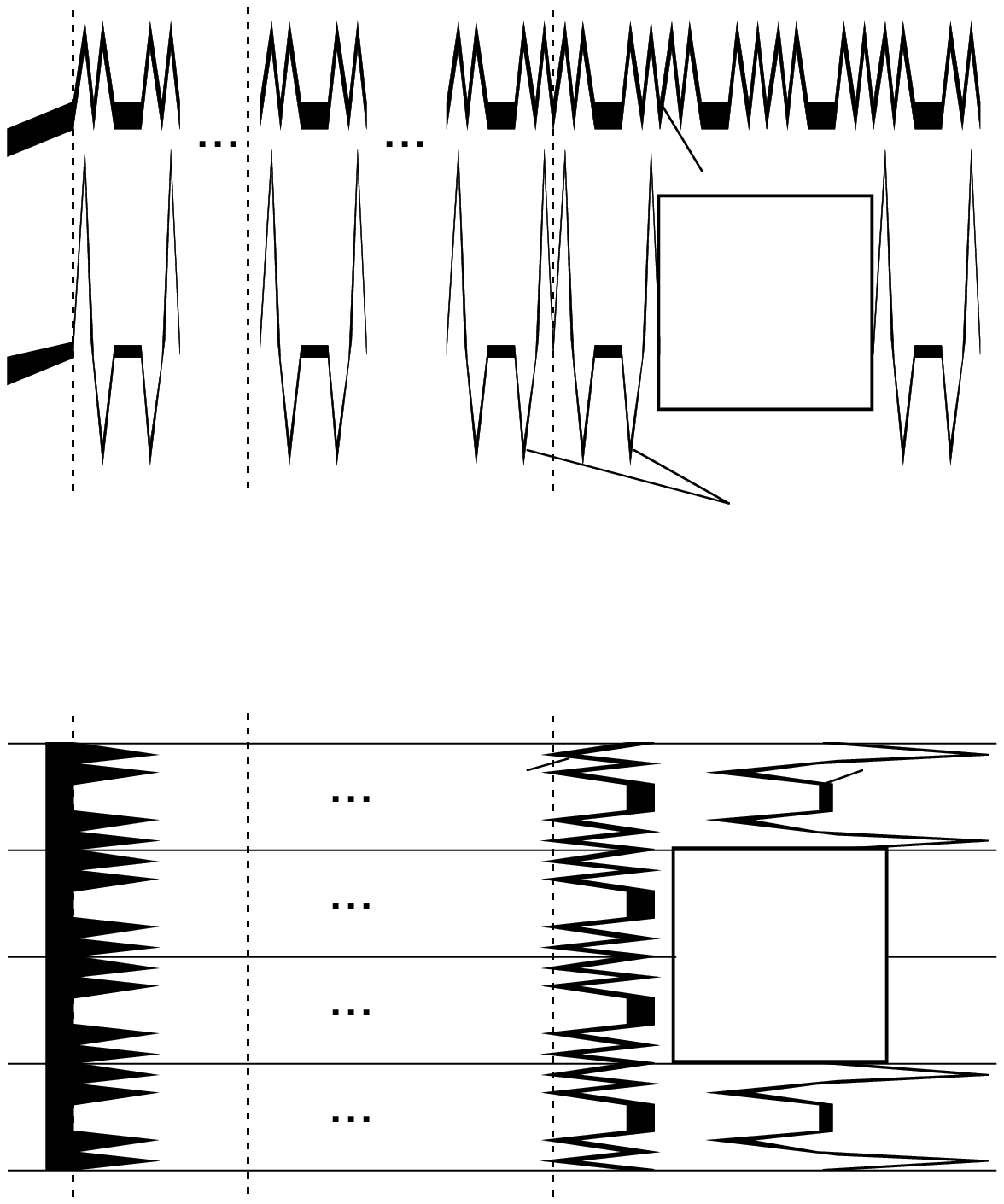}
\end{center}
\setlength{\unitlength}{.01in}
\hspace*{-4.98in}\begin{picture}(0,0)
\put(40,50){\makebox(0,0)[tl]{$t_{n}$}}
\put(120,50){\makebox(0,0)[tl]{$t_{n-1}$}}
\put(260,50){\makebox(0,0)[tl]{$t_{n-5}$}}
\put(180,360){\makebox(0,0)[tl]{\large Vertical Division}}
\put(180,22){\makebox(0,0)[bl]{\large Horizontal Division}}
\put(40,375){\makebox(0,0)[tl]{$t_{n}$}}
\put(120,375){\makebox(0,0)[tl]{$t_{n-1}$}}
\put(260,375){\makebox(0,0)[tl]{$t_{n-5}$}}
\put(345,530){\makebox(0,0)[tl]{\small Type 1 Cell}}
\put(360,380){\makebox(0,0)[tl]{\small Type 2 Cells}}
\put(190,250){\makebox(0,0)[tl]{\small Type 1 Cell}}
\put(420,250){\makebox(0,0)[tl]{\small Type 2 Cell}}
\end{picture}
\caption{Type~2 cells occur alongside other Type 2 cells or alongside
the holes $w_{n-1}\in W_{n-1}$
introduced in the
previous iteration.}
\label{fig:wheretps}
\end{figure}

\begin{figure}[bt]
\begin{center}
\includegraphics*{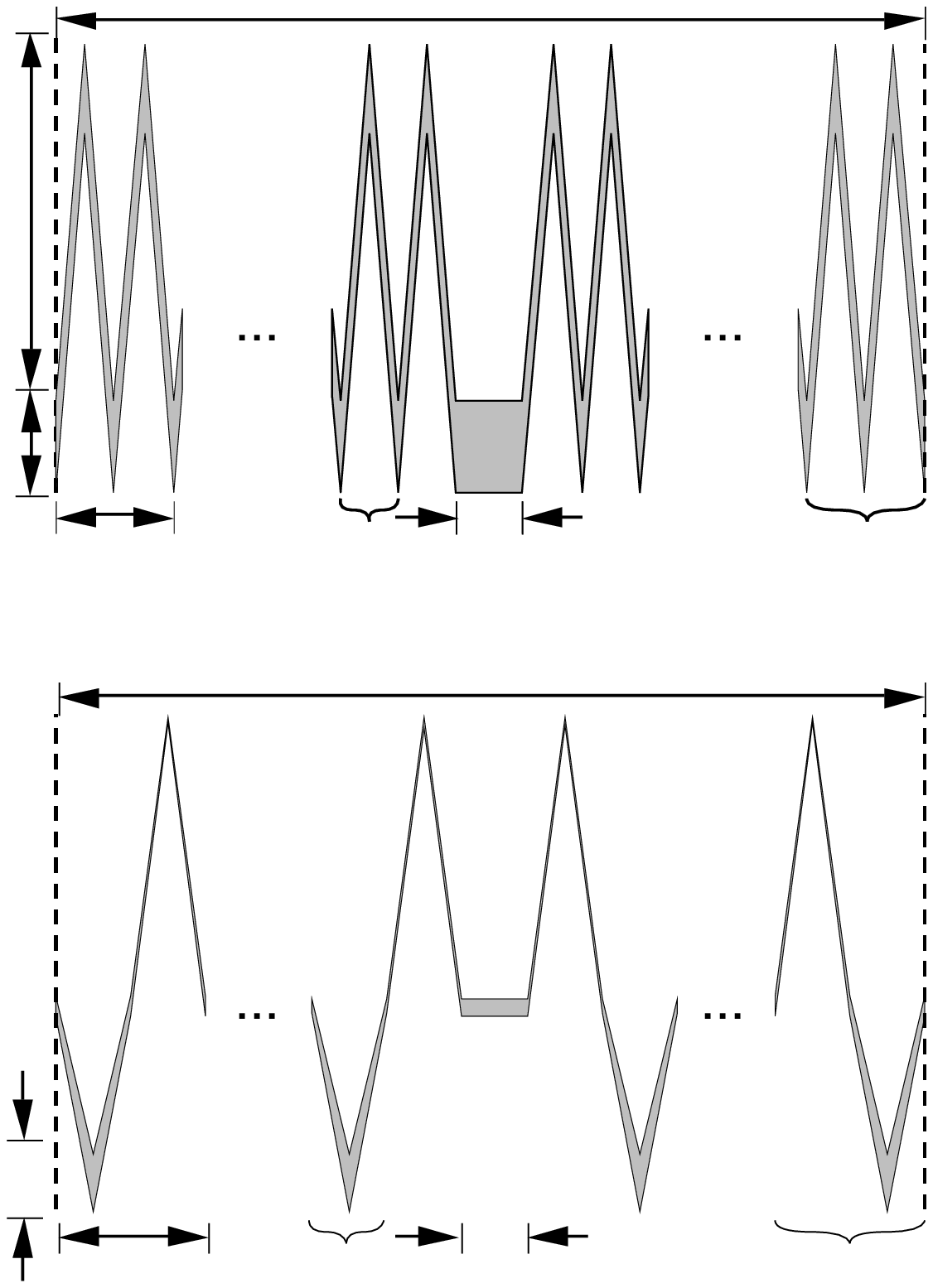}
\end{center}
\setlength{\unitlength}{.01in}
\hspace*{-4.81in}\begin{picture}(0,0)
\put(220,387){\makebox(0,0)[tl]{\large Cell Type 1}}
\put(65,407){\makebox(0,0)[tl]{$\frac{a_n-s_n}{k_n}$}}
\put(180,407){\makebox(0,0)[tl]{\small Cell-point}}
\put(257,422){\makebox(0,0)[tl]{$s_n$}}
\put(415,407){\makebox(0,0)[tl]{\small Cell-piece}}
\put(10,466){\makebox(0,0)[tl]{$<b_n$}}
\put(20,587){\makebox(0,0)[tl]{$c_n$}}
\put(260,657){\makebox(0,0)[tl]{$a_n$}}
\put(220,42){\makebox(0,0)[tl]{\large Cell Type 2}}
\put(75,67){\makebox(0,0)[tl]{$\frac{a_n-s_n}{k_n}$}}
\put(170,64){\makebox(0,0)[tl]{\small Cell-point}}
\put(260,77){\makebox(0,0)[tl]{$s_n$}}
\put(420,64){\makebox(0,0)[tl]{\small Cell-piece}}
\put(15,112){\makebox(0,0)[tl]{$<b_n$}}
\put(260,347){\makebox(0,0)[tl]{$a_n$}}
\end{picture}
\vspace*{.1in}
\caption{The two types of cells which occur to the
right of $\vertline(t_n)$ are Type~1 and Type~2 cells.}
\label{fig:ctypes}
\end{figure}

Type~1 cells are those that lie to the right of
$\vertline(t_n)$ that are not Type~2 cells.
In general Type~2 cells are those
that lie to the right of
$\vertline(t_n)$ and
that lie along
a hole $w_{n-1}\in W_{n-1}$.
To be more specific we must consider the vertical and horizontal
cases separately.
When $n$ is odd, a Type~2 cell is a cell that either lies
to the right of
$\vertline(t_{n-5})$ and also lies along the
vertical boundary of
a hole $w_{n-1}\in W_{n-1}$ or
lies between $\vertline(t_{n})$ and
$\vertline(t_{n-5})$  and also shares a vertical
boundary with another Type~2 cell.
When $n$ is even, a Type~2 cell is one that lies
to the right of
$\vertline(t_{n-5})$ and also lies along the
horizontal boundary of
a hole $w_{n-1}\in W_{n-1}$.
See Figure~\ref{fig:wheretps}.
For vertical cells both Type~1 or Type~2 cells
will consist of $k_n/2$ congruent
pieces on the left joined by a rectangle of width
$s_n$ to $k_n/2$ congruent pieces on the right.
See Figure~\ref{fig:ctypes}.
Such a cell will be symmetrical about a vertical line running through
the center of the rectangle.
Each of the $k_n$ pieces, which we call a {\it cell-piece},
will consist of two symmetrical
parts called {\it cell-points}.
Note that the width of a cell-piece is $(a_n-s_n)/k_n$.
In Type~1 cells the two cell-points of a cell-piece are congruent and
``point'' in the same direction.
For a typical cell $q_n\in Q_n$ of Type~1 see Figure~\ref{fig:ctypes}.
Note that $a_n$ defines the width of $q_n$.
The cell has a height of at least $c_n$ but less than $b_n+c_n$.
The thickness; i.e.,
vertical transverse thickness,
of the cell is limited by $b_n$.
Horizontal cells are similar to vertical cells except they
are rotated by ninety degrees.
 
For a typical vertical cell $q_n\in Q_n$ of Type~2
see Figure~\ref{fig:ctypes}.
Like Type~1 cells, $a_n$ defines the width of $q_n$ and the thickness
of the cell is limited by $b_n$.
In Type~2 cells, however, the two cell-points of a cell-piece will not
in general be congruent nor will they
``point''
in the same direction.
In addition the height of Type~2 cells can be greater than $c_n+b_n$
and in fact can have height greater than $2c_n+s_{n-1}$ where $s_{n-1}$
is the length of a side of the square holes in $W_{n-1}$.
When a Type~2 cell lies along a hole in $W_{n-1}$,
its cell-pieces are forced to extend at least $c_n/2$,
but no more
than $c_n+b_n$, beyond the hole.
Note that, when $n$ is odd (respectively even), cells to the right of
$\vertline(t_{n})$ that share a common vertical (respectively horizontal) boundary
are congruent.

\begin{figure}[tb]
\begin{center}
\includegraphics*{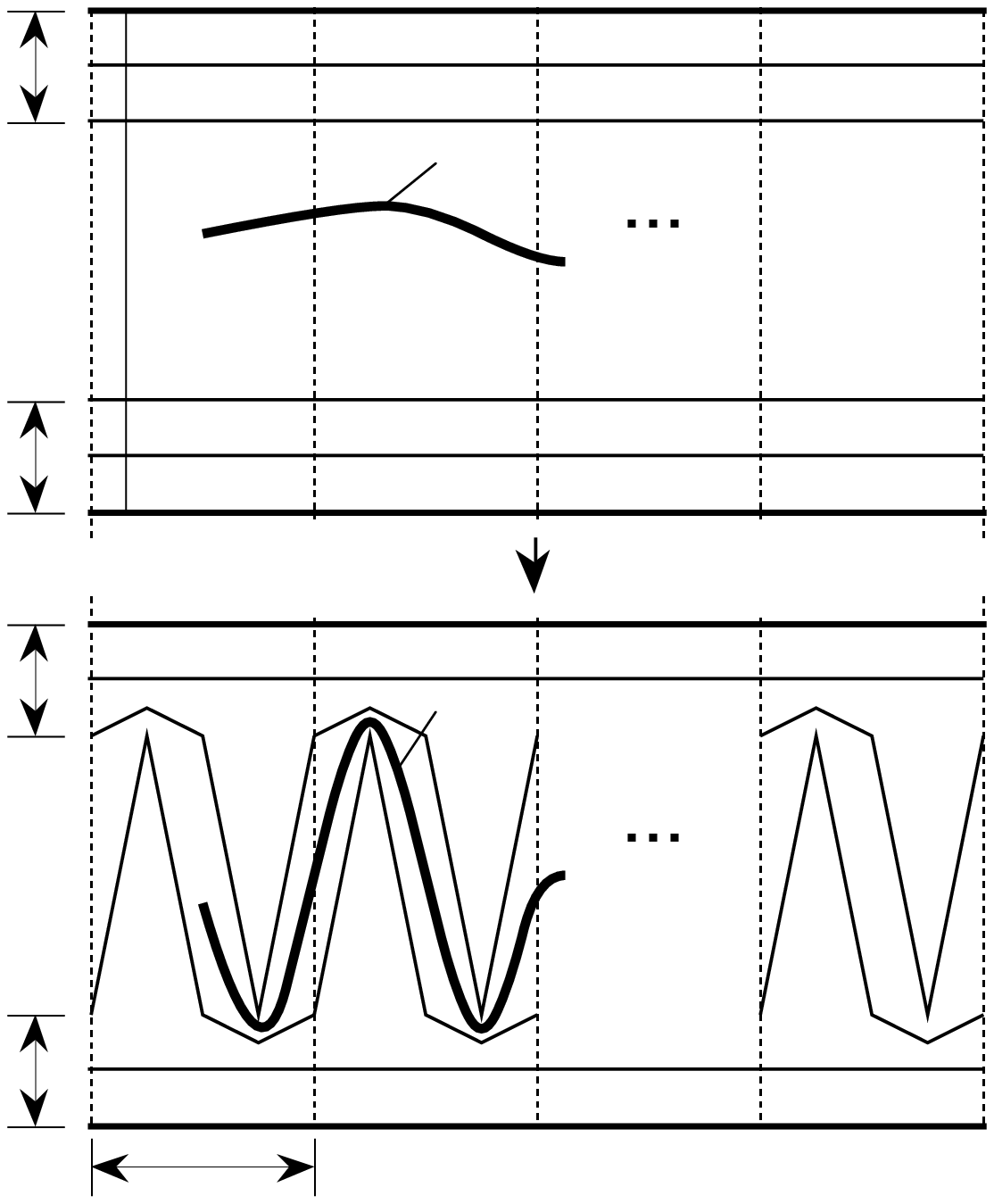}
\end{center}
\setlength{\unitlength}{.01in}
\hspace*{-4.91in}\begin{picture}(0,0)
\put(20,542){\makebox(0,0)[tl]{$a_n$}}
\put(225,497){\makebox(0,0)[tl]{$g$}}
\put(20,367){\makebox(0,0)[tl]{$a_n$}}
\put(53,322){\makebox(0,0)[tl]{$t_{n-3}$}}
\put(275,312){\makebox(0,0)[tl]{$f_n$}}
\put(453,322){\makebox(0,0)[tl]{$t_{n-4}$}}
\put(20,267){\makebox(0,0)[tl]{$a_n$}}
\put(230,257){\makebox(0,0)[tl]{$f_n(g)$}}
\put(20,92){\makebox(0,0)[tl]{$a_n$}}
\put(110,37){\makebox(0,0)[tl]{$a_n$}}
\end{picture}
\caption{The definition of $f_n$.}
\label{fig:deffn}
\end{figure}

\begin{figure}[tb]
\begin{center}
\includegraphics*{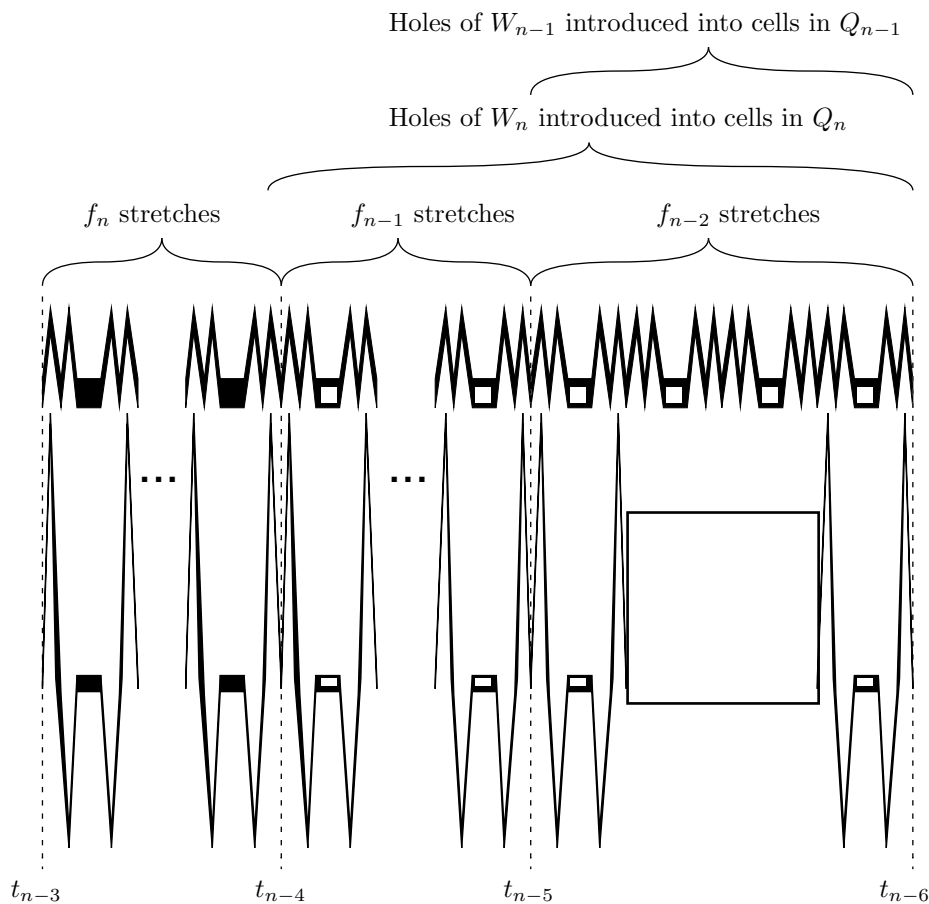}
\end{center}
\setlength{\unitlength}{.01in}
\hspace*{-5.1in}\begin{picture}(0,0)
\put(12,37){\makebox(0,0)[tl]{$t_{n-3}$}}
\put(140,37){\makebox(0,0)[tl]{$t_{n-4}$}}
\put(270,37){\makebox(0,0)[tl]{$t_{n-5}$}}
\put(467,37){\makebox(0,0)[tl]{$t_{n-6}$}}
\put(50,392){\makebox(0,0)[tl]{$f_{n}$ stretches}}
\put(190,392){\makebox(0,0)[tl]{$f_{n-1}$ stretches}}
\put(350,392){\makebox(0,0)[tl]{$f_{n-2}$ stretches}}
\put(210,442){\makebox(0,0)[tl]{Holes of $W_n$ introduced into cells in $Q_n$}}
\put(210,492){\makebox(0,0)[tl]{Holes of $W_{n-1}$ introduced into cells in $Q_{n-1}$}}
\end{picture}
\caption{Holes are only introduced into
cells in $Q_n$ to the right of $\vertline(t_{n-4})$.}
\label{fig:wherehls}
\end{figure}

At stage $n$ the function $f_n$ will be defined
so that eventually the elements of the 
decomposition that lie between $\vertline({t_{n-3}})$
and $\vertline({t_{n-4}})$
will be stretched vertically up close to the top edge of $D$ and
down towards the bottom edge of $D$.
Specifically, $f_n$ will
be defined so that a horizontal
line segment of length $a_n$
and lying between $\horiline(2-a_n)$ and
$\horiline(a_n)$ and between $\vertline({t_{n-3}})$
and $\vertline({t_{n-4}})$
will be bent above
$\horiline(2-a_n)$ and below $\horiline(a_n)$.
For $n>3$ the function
$f_n$ will be the identity to the right of
$\vertline({t_{n-4}})$ and to the left of $\vertline({t_{n-3}})$.
Between $\vertline(t_{n-3})$ and $\vertline(t_{n-4})$,
$f_n$ is a homeomorphism from
$[t_{n-3},t_{n-4}]\times [0,2]$ onto $[t_{n-3},t_{n-4}]\times [0,2]$
that maps vertical lines onto vertical lines and that is periodic in its first argument,
the period being $a_n$.
In addition $f_n$ between $\vertline(t_{n-3})$ and $\vertline(t_{n-4})$
satisfies the following constraints:
\begin{enumerate}
\item $f_n$ maps some point of the segment
$[t_{n-3},t_{n-3}+\frac{a_n}{2}]\times \{a_n\}$ above
$\horiline(2-a_n)$, and
\item $f_n$ maps some point of the segment
$[t_{n-3}+\frac{a_n}{2},t_{n-3}+a_n]\times \{2-a_n\}$
below $\horiline(a_n)$. 
\end{enumerate}
See Figure~\ref{fig:deffn}.
We define ${f_1=f_2=f_3=\identity_D}$.
Assuming that the sequence of functions $\{f_n\}_{n=1}^\infty$ has
been defined we define 
${F:(D\backslash E)\rightarrow(D\backslash E)}$
for every $x\in(D\backslash E)$ by the following equation
when $x$ lies between $\vertline(t_{n})$ and $\vertline(t_{n-1})$:
$$
F(x)=f_{n+3}\circ\cdot\cdot\cdot\circ f_1(x).
$$
Thus $F$ is a homeomorphism.
Once $f_n$ is defined we know how small the holes
must be to keep them from getting
stretched excessively by $f_n\circ f_{n-1}\circ\cdot\cdot\cdot\circ f_1$.
We will
only introduce these new holes to the right of $\vertline({t_{n-4}})$.
See Figure~\ref{fig:wherehls}.

\section{The Construction}
We now describe the construction in more detail.
We assume we are at stage $n$ where $n>1$ and that we are given
a division of $D$, $\widehat{R}_n$, and a continuum, $X_n\subset D$,
which is essentially a rectangle minus a finite number of open squares with
disjoint boundaries.
The smallest holes; i.e., those removed during the last stage, will be those
in $W_{n-1}$ and occur to the right of $\vertline(t_{n-5})$.
To start, we refine $\widehat{R}_n$ to create $R_n$ letting
$a_n=\mesh(\widehat{R}_n)/4$ be the mesh of $R_n$.
The size and position of the holes of $X_n$ will have been previously chosen
carefully so that the edges of the holes will lie on the
edges of the strips of $R_n$.

\subsection{The Creation of $Q_n$}\label{subsec:creation}

Our strategy in defining the cells of $Q_n$ will be to define a sequence
of disjoint polygonal lines $\{L^j_n\}_{j=-1}^{m_n+1}$
which run horizontally across the
vertical strips of $R_n$ when $n$ is odd and
run vertically across the
horizontal strips of $R_n$ when $n$ is even.
When $L^j_n$ runs horizontally, then by
$L^j_n(x)$ we will mean the $y$ such that $(x,y)\in L^j_n$.
Similarly, when $L^j_n$ runs vertically, then $(L^j_n(y),y)\in L^j_n$.
When $n$ is odd
we will define the
open cell $\widehat{q}_{i,j}$
for each $i\in1,...,{\frac{1}{2a_n}}$ and
$j\in0,...,(m_n+1)$
as follows:
$$
\widehat{q}_{i,j}=\{(x,y)\in D:(i-1)a_n< (x-\frac{1}{2})< ia_n
\mbox{ and }
L_n^{j-1}(x)< y< L_n^j(x)\}.
$$
We define $Q_n$ to be the closure of the nonempty cells:
$$
Q_n = \{\closure(\widehat{q}_{i,j}):
X_n\cap\widehat{q}_{i,j}\ne\emptyset
\mbox{, }
i\in1,...,\frac{1}{2a_n}
\mbox{ and }
j\in0,...,(m_n+1)
\}.
$$
When $n$ is even we define $Q_n$ in an analogous manner.
Because our strategy
for defining $\{L^j_n\}_{j=-1}^{m_n+1}$
is slightly different when $n$ is odd from when
$n$ is even we discuss them separately.

\subsubsection{$n$ is odd.}  The lines $\{L^j_n\}_{j=-1}^{m_n+1}$ will
be defined first between $\vertline(t_n)$ and
$\vertline(1)$, then between $\vertline(1/2)$ and
$\vertline(t_{n+1})$, and finally
between $\vertline(t_{n+1})$ and
$\vertline(t_n)$.  Between 
$\vertline(t_n)$ and
$\vertline(1)$ we proceed as in
\cite{Seaquist1995b}.
Let $O_n$ be the ordinates
of the top and bottom edges of the holes in $X_n$ along with $0$ and $2$.
For $y\in O_n\backslash\{2\}$
we denote by $\successor(y)$ the least element of $O_n$
greater than $y$.
Based on the same patterns as used in \cite{Seaquist1995b},
we define for each $y\in O_n$ three polygonal lines:  $\underline{M}^y_n$,
$M^y_n$, $\overline{M}^y_n$.
For each $y\in O_n\backslash\{2\}$ additional polygonal lines
are added between $M^y_n$ and $M^{\successor(y)}_n$ so that
\begin{enumerate}
\item the resulting cells are like the Type~1 and Type~2
cells described above;
\item  the number of rows of cells between 
$M^y_n$ and $M^{\successor(y)}_n$, which we
will denote by $m_y$, is such that
$(\successor(y)-y)/m_y$ is a constant,
which we denote by $d_n$; and
\item the constant $d_n$ divides
$2^{-(n+9)}$.
\end{enumerate}
Details on how this is done are contained in \cite{Seaquist1995b}.
We will in addition assume that $L^0_n(x)=\overline{M}^0_n(x)-b_n+d_n$ and
that $L^{m_n}_n(x)=\underline{M}^2_n(x)+b_n-d_n$.
Since $d_n<b_n$ this can be done by simply
inserting these two additional
polygonal lines where they are needed.
Note that none of these polygonal lines intersect and so there is
a natural order induced by where they intersect $\vertline(1)$.
Let $m_n=2/d_n-2$.
Reindex the lines in ascending order starting
with $L^{-1}_n=0$ and ending with $L^{m_n+1}_n=2$.

For each $x\in[1/2,t_{n+1}]$ and
for each $j\in\{-1,0,1,...,m_n+1\}$ we define
$L^j_n(x)=(j+1)d_n$.
For $x\in(t_{n+1},t_{n})$
we define
$$L^j_n(x)=
\frac{L^j_n(t_{n})-L^j_n(t_{n+1})}{t_{n}-t_{n+1}}
(x-t_{n+1})+L^j_n(t_{n+1}).$$

\subsubsection{$n$ is even.}  Again we define the polygonal
lines between
$\vertline(t_n)$ and $\vertline(1)$ first.
When $n>4$ we do this in two separate steps:  first between
$\vertline(t_{n-4})$ and $\vertline(1)$, and then between
$\vertline(t_n)$ and $\vertline(t_{n-4})$.
We
let $O_n$ be the
abscissa of all the vertical boundaries of the holes
in $X_n$ along with $t_{n-4}$
and $1$.
Based on the same pattern as used in \cite{Seaquist1995b}
we define for each $x\in O_n$
three polygonal lines:  $\underline{M}^x_n$,
$M^x_n$, $\overline{M}^x_n$.
For all $x\in O_n\backslash\{1\}$ additional polygonal lines
are added between $M^x_n$ and $M^{\successor(x)}_n$ so that
\begin{enumerate}
\item the resulting cells are like the Type~1 and Type~2
cells described above;
\item  the number of columns of cells between 
$M^x_n$ and $M^{\successor(x)}_n$, which we
will denote by $m_x$, is such that
$(\successor(x)-x)/m_x$ is a constant,
which we denote by $d_n$; and
\item the constant $d_n$ divides the following numbers 
$2^{-(n+9)}$, $(1-t_{n-4})$, $(t_{n-4}-t_{n})$, and $(t_n-1/2)$.
\end{enumerate}
Let $m'_n=((1-t_{n-4})/d_n)-2$
and re-index these polygonal lines between
$\vertline(1)$ and $\vertline(t_{n-4})$
starting on the right side with $L_n^{-1}=\vertline(1)$,
and continuing left to $L_n^{m'_n}=\overline{M}_n^{t_{n-4}}$ and
$L_n^{m'_n+1}=M_n^{t_{n-4}}$.
Note that $M_n^{t_{n-4}}\ne\vertline(t_{n-4})$,
that $M_n^{t_{n-4}}$ is just $\overline{M}^{t_{n-4}}_n$ shifted to the left
by $b_n$, and that we are indexing the lines in reverse order with
respect to that in \cite{Seaquist1995b}.
We continue defining the polygonal lines between
$\vertline(t_n)$ and $\vertline(t_{n-4})$ by
setting $L_n^{m'_n+1+i}(y)=L_m^{m'_n+1}(y)-i\cdot d_n$ for all
$i\in\{1,2,...,(t_{n-4}-t_n)/d_n\}$.
Set $m''_n=m'_n+(t_{n-4}-t_n)/d_n$.

When $n\le4$ let $t'=1-(a_{n-1}/2)$ and we again proceed in two separate steps:  first between
$\vertline(t')$ and $\vertline(1)$, and then between
$\vertline(t_n)$ and $\vertline(t')$.
We
let $O_n$ consist of $t'$
and $1$.
Based on the same pattern as used in \cite{Seaquist1995b}
we define for each $x\in O_n$ three polygonal lines:  $\underline{M}^x_n$,
$M^x_n$, $\overline{M}^x_n$.
For all $x\in O_n\backslash\{1\}$ additional polygonal lines
are added between $M^x_n$ and $M^{\successor(x)}_n$ so that
\begin{enumerate}
\item the resulting cells are like the Type~1
cells described above;
\item  the number of columns of cells between 
$M^x_n$ and $M^{\successor(x)}_n$, which we
will denote by $m_x$, is such that
$(\successor(x)-x)/m_x$ is a constant,
which we denote by $d_n$; and
\item the constant $d_n$ divides the following numbers:
$2^{-(n+9)}$, $a_{n-1}/2$, $(t'-t_n)$, and $(t_n-1/2)$.
\end{enumerate}
Let $m'_n=a_{n-1}/(2d_n)-2$
and re-index these polygonal lines between
$\vertline(1)$ and $\vertline(t')$
starting on the right side with $L_n^{-1}=\vertline(1)$,
and continuing left to $L_n^{m'_n}=\overline{M}_n^{t'}$ and
$L_n^{m'_n+1}=M_n^{t'}$.
We continue defining the polygonal lines between
$\vertline(t_n)$ and $\vertline(t')$ by
setting $L_n^{m'_n+1+i}(y)=L_m^{m'_n+1}(y)-i\cdot d_n$ for all
$i\in\{1,2,...,(t'-t_n)/d_n\}$.
Set $m''_n=m'_n+(t'-t_n)/d_n$.

Once we have defined all the polygonal lines
between 
$\vertline(t_n)$ and $\vertline(1)$ we define
the polygonal lines between 
$\vertline(1/2)$ and $\vertline(t_n)$ by
setting
$L_n^{m''_n+1+i}(y)=t_n-i\cdot d_n$
for all
$i\in\{1,2,...,(t_n-1/2)/d_n\}$.
Set $m_n=m''_n+((t_n-1/2)/d_n)$.

\subsection{Definition of the Homeomorphisms}\label{subsec:homeo} 
\par
We define $h_n:D\rightarrow D$ to be a
homeomorphism which maps vertical lines when $n$ is odd and horizontal
lines when $n$ is even onto themselves each
in a piecewise linear manner
so that
the preimage of the polygonal arcs $\{L^j_n\}_{j=-1}^{m_n+1}$ is
a collection of parallel straight lines evenly spaced apart at the
distance $d_n$.
Thus $h^{-1}_n(Q_n)$ is a collection of rectangles
with disjoint interiors.
Note that $h_n$ is the identity between $\vertline(1/2)$ and
$\vertline(t_{n+1})$.
Also, note that because of the way the
polygonal arcs $\{L^j_n\}_{j=-1}^{m_n+1}$ were defined,
$h_n$ maps
the boundary of each hole of $X_n$ onto itself.
Thus, if $w_i\in W_i$, then
$h_n(\closure(w_i))=\closure(w_i)$ for all $i<n$.

When $n$ is even there is another important result that follows
from the way the polygonal arcs $\{L^j_n\}_{j=-1}^{m_n+1}$
are defined; namely, that $\rwidth(x,h_n(x))<2a_{n-1}$.
This follows from the way the parameters
$a_n$, $b_n$, $c_n$ are chosen and from the following observations.
When $n>4$, then to the right of
$\vertline(t_{n-4})$
both $h^{-1}_n(L^i_n)$ and $L^i_n$ either intersect the same column of holes
in $W_{n-1}$ or they are both between two columns of holes from $W_{n-1}$
that are spaced no further than $a_{n-1}$ apart.
If
$h_n^{-1}(L_n^i)$ is to the left
of $\vertline(t_{n-4})$, then $h_n^{-1}(L_n^i)\cap L_n^i\ne\emptyset$.
Similarly when $n\le 4$, either $h_n^{-1}(L_n^i)\cap L_n^i\ne\emptyset$
or both $h_n^{-1}(L_n^i)$and $L_n^i$ are between $\vertline(1-a_{n-1}/2)$ and
$\vertline(1)$.

As in the constructions described in \cite{Seaquist1995a,Seaquist1995b}
the map $h_n$ can cause a great deal of stretching of cells that
lie along the horizontal edge of a hole (or $D$) when $n$ is odd or
along a vertical edge of a hole (or $D$) when $n$ is even.
If we are not careful, 
this stretching could potentially cause holes
introduced at stage $n$
to be overly enlarged
or to be too widely separated.
In an approach analogous to that described
in \cite{Seaquist1995a,Seaquist1995b} we
control where this stretching can occur in order to avoid problems.
We force it to occur at a distance between $a'_n/128$ and
$a'_n/64$ from the horizontal (vertical) boundary of $X_n$ when $n$ is odd (even).
When $n$ is odd $h_n$ is the identity below than $\horiline(a'_n/128)$
and above $\horiline(2-a'_n/128)$.
There is another place where $h_n$ can cause
a great deal of stretching and that is among the transition cells.
Note that in this case holes are only added to the transition
cells at a much later stage and that then the exact amount
of stretching is known and the holes can be made appropriately small
and close together.
We define $H_n$ to be $h_1\circ\cdot\cdot\cdot\circ h_{n-1}$ when
$n>1$ and $H_1$ to be $\identity_D$. 
Thus $Y_1=H_1(X_1)=D$
and $P_1=H_1(Q_1)=Q_1$.

\subsection{Preparation of the Next Stage} 
\par
We will set 
$$
P_n=H_n(Q_n)
\mbox{\quad and \quad}
Y_n=H_n(X_n).
$$
To continue the construction, we set
$$
W'_n=\{w_n:\exists q_n \in Q_n
\mbox{ and }
w_n
\mbox{ is an open }
s_n
\mbox{ by }
s_n\mbox{ square centered in }
h_n^{-1}(q_n)\}.
$$
We define $W_1=W_2=W_3=\emptyset$ and for $n>3$ we set
$$W_n=W'_n\cap([t_{n-4},1]\times[0,2]).$$
Define
$$X_{n+1}=X_n\backslash W^*_n.$$
This definition of $X_{n+1}$ prevents us from removing any holes to
the left of $\vertline(t_{n-4})$.
The next division $\widehat{R}_{n+1}$ is derived from projecting the
vertices of the rectangles $h^{-1}_n(Q_n)$ onto the $y$-axis if $n$ is odd
and onto the $x$-axis if $n$ is even.
Note that $\mesh(\widehat{R}_{n+1})=d_n$.

\section{The Specific Construction}
\par
We will now apply the construction to build the continuum
$Y=\cap_{n=1}^\infty Y_n$
and the collection $G=\{\cap_{n=1}^\infty \sstar(p_n,P_n)^*:
\cap_{n=1}^\infty p_n \ne\emptyset\mbox{ where } p_n\in P_n\}$ of subsets
of $Y$.
We first describe exactly how the parameters that control the construction
are chosen.  Assuming we are at stage $n$ the parameters are chosen in the
order given below.  The sequences $\{L_i\}_{i=-1}^\infty$ and
$\{K_i\}_{i=-1}^\infty$ are used to control the size of the
elements of the decomposition $G$.

\begin{enumerate}
\item Let $t_n={1/2}+(1/2)^{n+1}$ with $t_0=1$.
\item Let $a_n=a'_{n-1}$ with $a_1=(1/2)^{33}$.
\item Let $c_n=a_{n-1}/9$ with $c_1=(1/2)^{25}$.
\item Let $L_n=a_{n}/24$ with $L_{0}=L_{-1}=1/128$.
\item Let $K_n=4L_n$.
\item Define $f_n$ as described above so that between $\vertline(t_{n-3})$
and $\vertline(t_{n-4})$
a horizontal line segment that lies
between
$\horiline(a_n)$
and 
$\horiline(2-a_n)$
that is of width $a_n$ gets
bent above $\horiline(2-a_n)$ and
below $\horiline(a_n)$.
Note
that $f_1=f_2=f_3=\identity_D$.
As above define $F:(Y_1\backslash E)\rightarrow (Y_1\backslash E)$
by $F(x)=f_{n}\circ\cdot\cdot\cdot\circ f_{1}(x)$ if $x$ is between
$\vertline(t_{n-3})$ and $\vertline(t_{n-4})$
for $n>3$.  
\item Let $\delta_n>0$ so that
\begin{enumerate}
\item $|x-x'|<\delta_n\implies |H_n(x)-H_n(x')|<L_n/2^{n+1}$, and
\item $|x-x'|<\delta_n\implies |
f_n\circ\cdot\cdot\cdot\circ f_{1}\circ H_n(x)-
f_n\circ\cdot\cdot\cdot\circ f_{1}\circ H_n(x')|<(1/2)^{n+1}$.
\end{enumerate}
\item Let $b_n>0$ be rational so that
\begin{enumerate}
\item  $b_n < \delta_n/4$,
\item  $b_n < a_{n}/(2k_{n-1})$ where $k_{0}=512$,
\item  $b_n < b_{n-1}(a_{n-1}-s_{n-1})/(4k_{n-1}(c_{n-1}+b_{n-1}))$
when $n>1$, and
\item  $b_n < a_{n}/2^{n+6}$.
\end{enumerate}
\item Define $m_n$ and $d_n$ as in~\ref{subsec:creation}.  Recall that we will have
\begin{enumerate}
\item  $d_n\mbox{ divides }(1/2)^{n+9}$ and
$d_n\mbox{ divides }a_{n-1}/2$.
\end{enumerate}
\item Define $h_n$ as in~\ref{subsec:homeo}.
Define $H_n=h_1\circ\cdot\cdot\cdot\circ h_{n-1}$ with $H_1=\identity_D$.
\item Let $a'_n=d_n/4$ and $s_n=2a'_n=d_n/2$.
\item Let $k_n$ be an integer so that
\begin{enumerate}
\item $k_n\ge512$,
\item $4\mbox{ divides }k_n$, and
\item $k_n>a_n/a'_n$.
\end{enumerate}
\end{enumerate}

It can be shown as in \cite{Seaquist1995b} that at each stage $n$
we can pick $a_n$, $a'_n$, $b_n$, $c_n$, $s_n$, and $k_n$ following the
above constraints so that
it is possible to continue the construction to stage $n+1$.
We must also be able to show that
\begin{enumerate}
\item  the continuum $Y=\cap_{n=1}^\infty Y_n$ is the Sierpi\'nski
curve;
\item  the continuum $F(Y\backslash E)\cup E$ is the Sierpi\'nski curve;
\item  all members of $G$  except those that intersect the left edge of $D$
are nondegenerate;
\item  the collection $G=\{\cap_{n=1}^\infty \sstar(p_n,P_n)^*:
\cap_{n=1}^\infty p_n \ne\emptyset\mbox{ where } p_n\in P_n\}$ is a continuous
decomposition of $Y$. 
\item  the collection
$G'=\{F(g)|g\in G \mbox{ and } g\cap E=\emptyset\}\cup\{E\}$
is a continuous decomposition of $F(Y\backslash E)\cup E$; and
\item the decomposition space $(F(Y\backslash E)\cup E)/G'$ is homeomorphic to
the Sierpi\'nski curve.
\end{enumerate}
We can show that 1--4 above hold in precisely the same way as they are shown
in~\cite{Seaquist1995b}.  To show 5--6 we need the following lemma. 
It guarantees that every $g\in G$ will be bent and
sufficiently stretched vertically by $F$.

\begin{lemma}\label{lm:bentg}
Let $g\in G$ and let $n$ be the least integer such that $g$ is strictly
to the right of $\vertline(t_{n+1})$.
Let $g=\cap_{i=1}^\infty\sstar(p_i,P_i)^*$
where $p_i\in P_i$ for all $i\in\Z^+$ and $\cap_{i=1}^\infty p_i\ne\emptyset$.
Let $q_i=H^{-1}_i(p_i)$ for all $i$.
Assume $n>4$.
 
If $p_{n-1}\cap \vertline(t_{n+1})=\emptyset$,
then there exists a trapezoid $T$ between
$\vertline(t_{n+1})$ and $\vertline(t_{n-1})$
so that
\begin{enumerate}
\item [$1$.] $\rwidth(T)>a_{n+3}$,
\item [$2$.] if $(x,y)\in T$, then $a_{n+3}<y<2-a_{n+3}$,
\item [$3$.] every vertical line that intersects $T$ will intersect
$g\cap T$;
\end{enumerate}
and

if
$p_{n-1}\cap \vertline(t_{n+1})\ne\emptyset$,
then there exists a trapezoid $T$ between
$\vertline(t_{n+1})$ and $\vertline(t_{n})$
so that
\begin{enumerate}
\item [$1'$.] $\rwidth(T)>a_{n+4}$,
\item [$2'$.] if $(x,y)\in T$, then $a_{n+4}<y<2-a_{n+4}$,
\item [$3'$.] every vertical line that intersects $T$ will intersect
$g\cap T$;
\end{enumerate}
\end{lemma}

\begin{proof}

Since $H_{n-1}=\identity_D$ to the left of $\vertline(t_{n-1})$ it can be shown
that $p_{n-1}=q_{n-1}$
because of the
definition of $n$.
Since $p_{n-1}$ must be between $\vertline(t_{n+1})$ and $\vertline(t_{n-1})$, we know that
for any $m\in\Z^+$ there some point of $p_{n+m}$
that lies between
$\vertline(t_{n+1})$
and $\vertline(t_{n-1})$.
Thus for any $m\in\Z^+$ we have that
$q_{n+m}$ is strictly
to the right of $\vertline(t_{n+2})$.
(See, for example, the proof of Lemma~9 of~\cite{Seaquist1998}). 
Applying Claim~C of~\cite{Seaquist1998}
we have that for all $m\ge n$
\begin{eqnarray}\label{eqna:C}
p_{m+3}&\subset&\sstar(p_{m+4-1},P_{m+4-1})^*\\
&\subset&\sideset{}{_{\frac{12L_{m+2}+3K_{m+2}}{2^{m+4}}}}\neighborhood(g)\notag\\
&\subset&\sideset{}{_{\frac{a_{m+2}}{2^{m+4}}}}\neighborhood(g).\notag
\end{eqnarray}

First we assume that $p_{n-1}\cap\vertline(t_{n+1})=\emptyset$.
If $g$ lies between $\horiline(a_{n+3})$ and $\horiline(2-a_{n+3})$, then
we are done because it can be shown that $\rwidth(g)>a_{n+3}$.
See, for example, the proof of Lemma~10 in~\cite{Seaquist1998}.
So assume that there is a point $x\in g$ that is below $\horiline(a_{n+3})$.
The situation when
there is a point $x\in g$ that is above $\horiline(2-a_{n+3})$ is handled
similarly.

{\bf Case $1$:} Assume that $n$ is odd in addition to assuming $p_{n-1}\cap\vertline(t_{n+1})=\emptyset$.
Since $g\subset\sstar(p_{n+1},P_{n+1})^*$ there is 
a $\widehat{p}_{n+1}\in\sstar(p_{n+1},P_{n+1})$
so that $x\in\widehat{p}_{n+1}$.
Note that
$\widehat{q}_{n+1}=H^{-1}_{n+1}(\widehat{p}_{n+1})$ is a
horizontal cell since $n+1$ is even.
Because of the way polygonal lines were defined, specifically
that 
$L^0_n(x)=\overline{M}^0_n(x)-b_n+d_n$ and
that $L^{m_n}(x)_n=\underline{M}^2_n+b_n-d_n$, we know that $h_n$ is the identity
below $\horiline(d_n)$ and above $\horiline(2-d_n)$
between $\vertline(t_{n+1})$ and $\vertline(t_{n})$ and
because of the way stretching is controlled near boundaries
we know that
$h_{n}$ is the identity below
$\horiline(a_{n+1}/128)$.
Note that $h_{n-1}$
is the identity to the left of $\vertline(t_{n}+3d_{n-1})$.
Thus since $x$ is below $horiline(a_{n+3})$, we have that
$H^{-1}_{n+1}(x)$ and that $x\in\widehat{q}_{n+1}$.
Let $\widehat{q}^p_{n+1}$ be a cell-piece of the cell 
$\widehat{q}_{n+1}$ that lies between
$\horiline(a_{n+1}/128)$ and $\horiline(a_{n+1}/256)$.
Such a cell-piece exists
since $k_{n+1}\ge512$ and the width of a cell-piece is $(a_{n+1}-s_{n+1})/k_{n+1}$.
Let $\ell_1$ and $\ell_2$ be the horizontal lines that run
along the top and bottom boundaries of $\widehat{q}^p_{n+1}$
respectively.
Define $T'$ to be the part of $\sstar(q_{n+1},Q_{n+1})^*$
that lies between $\ell_1$ and $\ell_2$.
Notice that since $a_{n+3}<a_{n+1}/256$, $x$ is below $\ell_2$.
See Figure~\ref{fig:Tbox}.
\begin{figure}[tb]
\begin{center}
\includegraphics*{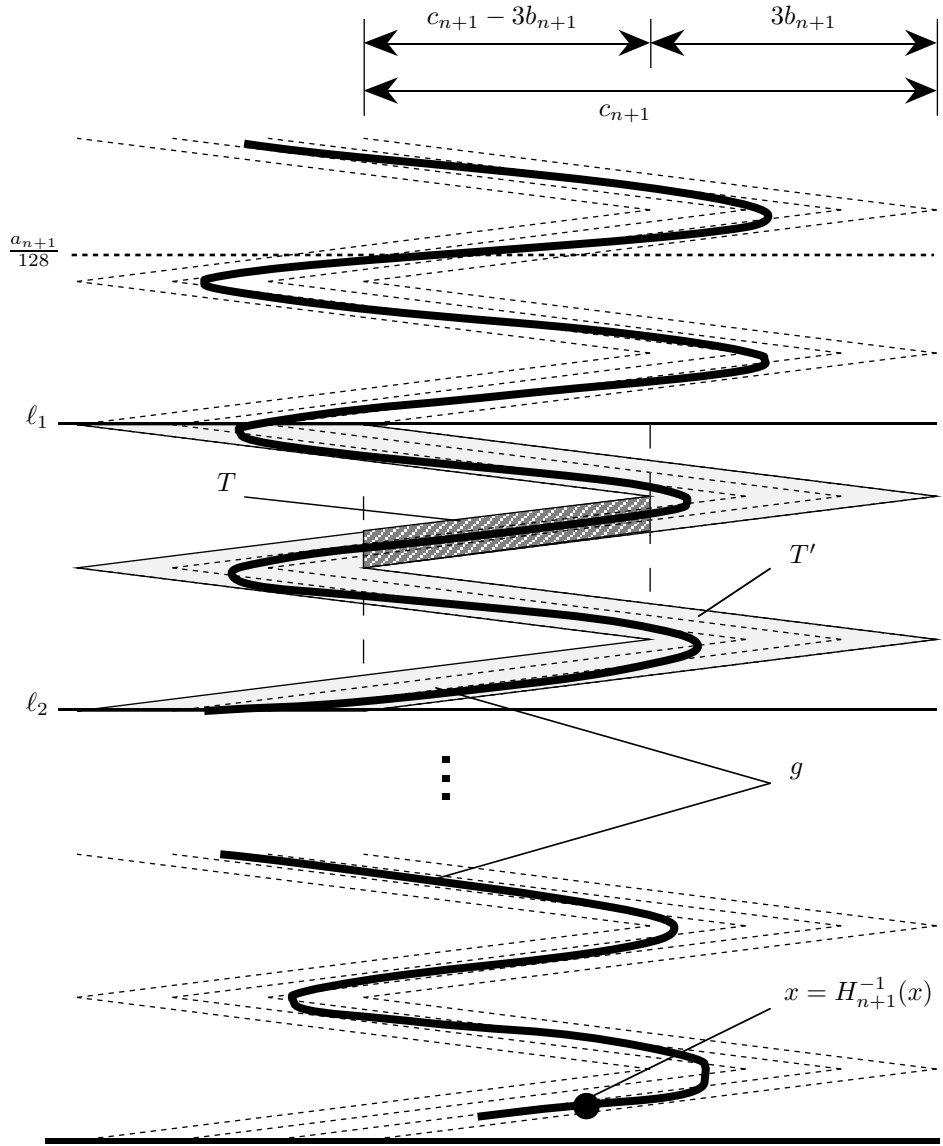}
\end{center}
\setlength{\unitlength}{.01in}
\hspace*{-5.1in}\begin{picture}(0,0)
\put(230,624){\makebox(0,0)[tl]{$c_{n+1}-3b_{n+1}$}}
\put(410,624){\makebox(0,0)[tl]{$3b_{n+1}$}}
\put(320,572){\makebox(0,0)[tl]{$c_{n+1}$}}
\put(10,503){\makebox(0,0)[tl]{$\frac{a_{n+1}}{128}$}}
\put(20,414){\makebox(0,0)[tl]{$\ell_1$}}
\put(120,379){\makebox(0,0)[tl]{$T$}}
\put(420,342){\makebox(0,0)[tl]{$T'$}}
\put(20,264){\makebox(0,0)[tl]{$\ell_2$}}
\put(420,227){\makebox(0,0)[tl]{$g$}}
\put(417,115){\makebox(0,0)[tl]{$x=H^{-1}_{n+1}(x)$}}
\end{picture}
\caption{Case 1 of Lemma~\ref{lm:bentg}.}
\label{fig:Tbox}
\end{figure}
Now $H_{n+1}(T')=H_{n-1}\circ h_{n-1}\circ h_{n}(T')=T'$,
and $T'\subset\sstar(p_{n+1},P_{n+1})^*$.

Consider $p_{n+3}$.
There is $q'_{n+2}\in\sstar(q_{n+2},Q_{n+2})$
so that $q_{n+3}$ crosses some cell-point of $q'_{n+2}$.
But $c_{n+2}=a_{n+1}/9>a_{n+1}/64$ so some point
of $p_{n+3}$ lies above
$\horiline(a_{n+1}/64)$ and since $a_{n+2}/2^{n+2}<a_{n+1}/128$ we
can conclude
from Equation~\ref{eqna:C} that some
point of $g$ lies above $\ell_1$.
Therefore $g$ must
lie along $T'$ and ``cross'' one of the ``cell points'' of $T'$.
Now $c_{n+1}-3b_{n+2}>a_{n+3}$ so there exists a trapezoid
$T\subset T'$ that satisfies conditions {1--3}.

{\bf Case $2$:}  Assume $n$ is even in addition to assuming $p_{n-1}\cap\vertline(t_{n+1})=\emptyset$.
Consider
$\widehat{p}_{n+2}\in\sstar(p_{n+2},P_{n+2})$
so that $x\in\widehat{p}_{n+2}$.
Note that
$\widehat{q}_{n+2}=H^{-1}_{n+2}(\widehat{p}_{n+2})$ is a
horizontal cell.
If $q'_{n+2}\in\sstar(q_{n+2},Q_{n+2})$ then denote $H_{n+2}(q'_{n+2})$ by $p'_{n+2}$.
Since $H_{n-2}$ is the identity to the left of $\vertline(t_{n-2})$
we have that $p'_{n+2}=h_{n-1}\circ h_{n}\circ h_{n+1}(q'_n+2)$.
Now $h_{n+1}$ is the identity below $\horiline(a_{n+2}/128)$.
Also we know that $h_{n-1}$ is the identity below $\horiline(d_{n-1})$ and
to the left of $\vertline(t_{n-1})$.
Recall that $d_{n-1}>a_{n+1}/128$.
Finally, because $h_{n}$ maps horizontal lines onto horizontal lines
and because the image under $h_{n}$ of any point of
$\sstar(q_{n+2},Q_{n+2})$ is to the left of $\vertline(t_{n-1})$,
we can conclude that $H_{n+2}(q'_{n+2})=h_n(q'_{n+2})$ and
that $h^{-1}_n(x)\in\widehat{q}_{n+2}$ which is below
$a_{n+3}$.
Let $\widehat{q}^p_{n+2}$ be a cell-piece of the cell
$\widehat{q}_{n+2}$ that lies between
$\horiline(a_{n+2}/128)$ and $\horiline(a_{n+2}/256)$.
Let $\ell_1$ and $\ell_2$ be the horizontal lines that run
along the top and bottom boundaries of $\widehat{q}^p_{n+2}$
respectively.
Define $T'$ to be the part of $\sstar(q_{n+2},Q_{n+2})^*$
that lies between $\ell_1$ and $\ell_2$.
Notice that since $a_{n+3}<a_{n+2}/256$, $x$ is below $\ell_2$
and that $H_{n+2}(T')=h_{n}(T')$.
Again from Equation~\ref{eqna:C}  and the fact that $c_{n+3}>a_{n+2}/64$
we can conclude that some
point of $g$ lies above $\ell_1$ and that $g$ is forced to lie
along the image under $h_n$ of a ``cell-piece'' of $T'$. 
Because of the behavior of $h_{n}$ to the left of $\vertline(t_{n-1})$ we can find
a trapezoid $T\subset h_{n}(T')$ that satisfies conditions
{1--3} of the lemma.

Now we will assume that $p_{n-1}\cap\vertline(t_{n+1})\ne\emptyset$.
We also assume that there is a point $x\in g$ that is below
$\horiline(a_{n+4})$.
The situation when $x$ is close to the top boundary of $D$ is handled
similarly.

{\bf Case $1'$:}  Assume that $n$ is odd in addition to assuming $p_{n-1}\cap\vertline(t_{n+1})\ne\emptyset$.
Consider
$\widehat{p}_{n+3}\in\sstar(p_{n+3},P_{n+3})$
so that $x\in\widehat{p}_{n+3}$.
Note that
$\widehat{q}_{n+3}=H^{-1}_{n+3}(\widehat{p}_{n+3})$ is a
horizontal cell.  
If $q'_{n+3}\in\sstar(q_{n+3},Q_{n+3})$ then denote $H_{n+3}(q'_{n+3})$ by $p'_{n+3}$.
Since $H_{n}$ is the identity to the left of $\vertline(t_{n})$
we have that $p'_{n+3}=h_{n}\circ h_{n+1}\circ h_{n+2}(q'_n+3)$.
As in Case~1 above we see that $h_{n+2}$ is the identity below $\horiline(a_{n+3}/128)$.
Also we know that $h_n$ is the identity below $\horiline(d_n)$ and to the
left of $\vertline(t_n)$.
Finally because $h_{n+1}$ maps horizontal lines onto horizontal lines
and because the image under $h_{n+1}$ of any point of
$h_{n+2}(\sstar(q_{n+3},Q_{n+3}))$ is to the left of $\vertline(t_n)$
we can conclude that $H_{n+3}(q'_{n+3})=h_{n+1}(q'_{n+3})$ and
that $h^{-1}_{n+1}(x)\in\widehat{q}_{n+3}$ which is below
$a_{n+4}$.
Let $\widehat{q}^p_{n+3}$ be a cell-piece of the cell 
$\widehat{q}_{n+3}$ that lies between
$\horiline(a_{n+3}/128)$ and $\horiline(a_{n+3}/256)$.
Let $\ell_1$ and $\ell_2$ be the horizontal lines that run
along the top and bottom boundaries of $\widehat{q}^p_{n+3}$
respectively.
Define $T'$ to be the part of $\sstar(q_{n+3},Q_{n+3})^*$
that lies between $\ell_1$ and $\ell_2$.
Notice that since $a_{n+4}<a_{n+3}/256$, $x$ is below $\ell_2$, and
that $H_{n+3}(T')=h_{n+1}(T')$.
From Equation~\ref{eqna:C} and the fact that $c_{n+4}>a_{n+3}/64$
we can conclude that some
point of $g$ lies above $\ell_1$ and that $g$ is forced to lie
along the image under $h_{n+1}$ of a ``cell-piece'' of $T'$. 
Because of the behavior of $h_{n+1}$ near $\vertline(t_{n+1})$ we can find
a rectangle $T\subset h_{n+1}(T')$ that satisfies conditions
{$1'$--$3'$} of the lemma.

{\bf Case $2'$:}   Assume $n$ is odd in addition to assuming $p_{n-1}\cap\vertline(t_{n+1})\ne\emptyset$.
Exactly as in Case~2 above.
\end{proof}

We now observe that because $G$ is a continuous decomposition
and the elements of $G$ that intersect the left edge of $D$
are singletons, the elements of $G$ have smaller and smaller diameters
as they get closer to $E$.
Because of the definition of $F$ for each $g\in G$ with
$g\cap E=\emptyset$ we have that $\rwidth(g)=\rwidth(F(g))$.
Thus the elements $g'\in G'$ have a smaller
and smaller width as they get closer and closer to $E$.

\begin{lemma}\label{lm:Gpcont}
The collection
$G'=\{F(g):g\in G\mbox{ and } g\cap E=\emptyset\}\cup \{E\}$
is a continuous decomposition of $Y'=F(Y\backslash E)\cup E$.
\end{lemma}
\begin{proof}
Since we have that $G$ is a continuous decomposition of $Y$
and that $F|(Y\backslash E)$ is a homeomorphism, then $G'$ must be continuous
at every $g'\in G'$ except possibly when $g'=E$.  But
the Lemma~\ref{lm:bentg} guarantees that the elements of $G'$ get
bent along $E$ so $G'$ is lower semicontinuous
at $E$ and by the argument preceding this lemma  we have that
$G'$ is upper semicontinuous at $E$.  Thus $G'$ is continuous
at $E$ also. 
\end{proof}

\begin{lemma}\label{lm:Gpsier}
The set $Y'/G'$ is the Sierpi\'nski curve.
\end{lemma}

\begin{proof}
Let $\pi:Y'\rightarrow Y'/G'$
be the natural projection.
We know that $Y'/G'$ is locally connected because $Y'$ is.
Since $G'$ is upper semicontinuous, we know that by adding the points
of the compliment of $Y'$,
we get an upper semicontinuous decomposition
of ${\E}^2$ which we will call $G''$.
Thus we extend $\pi$
to $\widehat{\pi}:{\E}^2\rightarrow {\E}^2/G''$.
No member of $G''$ separates
the plane so by R.L. Moore's Theorem~\cite{Daverman1986}
we know that ${\E}^2/G''$ is homeomorphic to the plane.
Thus $Y'/G'$ is planar.
By Corollary~13 in~\cite{Seaquist1998}, we get
that
$(\widehat{\pi})|\closure(U)$ is a homeomorphism for any $U$ a bounded
component
of $(Y')^c$.
Finally, if $V$ is the unbounded component of $(Y')^c$,
then $\widehat{\pi}(\boundary(V))$ is also a simple closed curve.
Thus the images under $\widehat{\pi}$
of the boundaries of the components of $(Y')^c$ are simple
closed curves, are dense in $Y'/G'$, and have diameters that go to zero.
Therefore
$Y'/G'$ is homeomorphic to the Sierpi\'nski
curve by \cite{Whyburn1958a}. 
\end{proof}

This proves Theorem~\ref{thm:construction}.

\bibliographystyle{amsalpha}

\end{document}